\documentclass[10pt]{article}
\usepackage{amssymb}
\usepackage{amsfonts}

\title{Countable locally $2$-arc-transitive bipartite graphs}
\author{Robert D. Gray, University of East Anglia, and } 
\date{J K Truss, University of Leeds$^1$.}
\begin{document}
\maketitle 
\newtheorem{lemma}{Lemma}[section]
\newtheorem{theorem}[lemma]{Theorem}
\newtheorem{corollary}[lemma]{Corollary}
\newtheorem{definition}[lemma]{Definition}
\newtheorem{prop}[lemma]{Proposition}
\newtheorem{claim}[lemma]{Claim}
\newtheorem{assumption}[lemma]{Assumption}
\newtheorem{remark}[lemma]{Remark}
\newtheorem{conjecture}[lemma]{Conjecture}
\setcounter{footnote}{1}\footnotetext{This work was supported by EPSRC grants EP/D04829/1 and EP/H00677X/1. The first author was also 
supported by EPSRC Postdoctoral Fellowship EP/E043194/1 held at the University of St Andrews, Scotland, and was also partially supported by FCT and 
FEDER, project POCTI-ISFL-1-143 of Centro de \'{A}lgebra da Universidade de Lisboa, and by the project PTDC/MAT/69514/2006. \\
Keywords bipartite graph, locally 2-arc transitive, Dedekind--Macneille completion, 2010 Mathematics subject classification: 05C20, 05C38.}

\newcounter{number}

\begin{abstract} 
We present an order-theoretic approach to the study of countably infinite locally 2-arc-transitive bipartite graphs. Our approach is motivated by 
techniques developed by Warren and others during the study of cycle-free partial orders. We give several new families of previously unknown countably 
infinite locally-2-arc-transitive graphs, each family containing continuum many members. These examples are obtained by 
gluing together copies of incidence graphs of semilinear spaces, satisfying a certain symmetry property, in a tree-like way. 
In one case we show how the classification problem 
for that family relates to the problem of determining a certain family of highly arc-transitive digraphs.
Numerous illustrative examples are given.   
\end{abstract}

\maketitle

\section{Introduction}

The study of finite graphs satisfying certain transitivity properties has a long and distinguished history with much of the motivation coming from the 
world of permutation groups. Various symmetry conditions have been considered, in various contexts, the more important among them including 
$k$-arc-transitivity (see \cite{Biggs} and Tutte's seminal work \cite{Tutte1, Tutte2}), distance-transitivity  (and distance-regularity) (see 
\cite{Brouwer}), and homogeneity (see \cite{Gardiner}).

The theory of infinite graphs satisfying transitivity conditions is in some aspects less developed, and necessarily so, since the powerful tools of 
finite group theory (like the classification of finite simple groups and the O'Nan--Scott Theorem for quasi-primitive groups) are not available to us. 
In this area of groups acting on infinite graphs (and other infinite relational structures) techniques and ideas from model theory play a key role. 
Arising from this, the strongest condition that has received attention, called \emph{homogeneity}, has led to a classification for graphs (see 
\cite{Gardiner,Lachlan}) and for other structures such as posets and digraphs (\cite{Schmerl, Cherlin}). However, for the other conditions mentioned 
above, even in instances where the picture for finite graphs is fairly complete, the situation for infinite graphs remains mysterious. 

The objects of study in this paper will be the countable infinite {\em $2$-arc-transitive graphs} and more generally {\em locally $2$-arc-transitive 
graphs}. Given a graph $\Gamma = (V\Gamma, E\Gamma)$ with vertex set $V\Gamma$ and edge set $E\Gamma$, an \emph{$s$-arc} in $\Gamma$ is an 
$(s+1)$-tuple $v_0, v_1, \ldots, v_s$ of vertices such that $v_i$ is adjacent to $v_{i+1}$ and $v_i \neq v_{i+2}$ for all $i$. The graph $\Gamma$ is 
\emph{locally $s$-arc-transitive} if, with $G = {\rm Aut}(\Gamma)$,  for each vertex $v$, the vertex stabilizer $G_v$ is transitive on the set of 
$s$-arcs starting at $v$. If $G$ is transitive on the set of all $s$-arcs in $\Gamma$ we say that $\Gamma$ is \emph{$s$-arc-transitive}. Being 
$s$-arc-transitive is equivalent to being simultaneously both vertex transitive and locally $s$-arc-transitive. Interest in (finite) $s$-arc-transitive 
graphs goes back to the fundamental work of Tutte \cite{Tutte1, Tutte2} showing that $s$-arc-transitive graphs of valency 3 satisfy $s \leq 5$. 
Later Weiss \cite{Weiss1} proved that if the valency is at least 3, then $s \leq 7$. His result relies on the classification of finite simple groups, 
as do many other results in the area; see for example \cite{Li}. We remark that a graph is $s$-arc-transitive or locally $s$-arc transitive if and 
only if all its connected components are (and for $s$-arc-transitivity are all isomorphic), so we generally restrict to the case of connected graphs 
without further comment.

Let $\Gamma$ be a connected $s$-arc-transitive graph, and $G = {\rm Aut}(\Gamma)$. If $G$ does not act transitively on vertices and all vertices have 
valency at least two then, since $G$ is edge-transitive, it follows that $\Gamma$ is bipartite and $G$ has two orbits on $\Gamma$, which are precisely 
the blocks of the bipartition of $\Gamma$. The locally $s$-arc-transitive graphs have received a lot of attention in the literature partly due to 
their links with areas of mathematics such as generalized $n$-gons, groups with a $(B,N)$-pair of rank two, Moufang polygons and Tutte's $m$-cages; 
see \cite{Weiss1} and \cite{Delgado} for more details. Important recent papers about local $s$-arc-transitivity include 
\cite{Giudici1, Giudici2, Giudici3, Giudici4}. In particular in \cite{Giudici1} a programme of study of locally $s$-arc-transitive graphs for which 
$G$ acts intransitively on vertices was initiated, with the O'Nan-Scott Theorem for quasi-primitive groups playing an important role in their 
reduction theorem. Interesting connections with semilinear spaces and homogeneous factorizations were explored in \cite{Giudici3}. 

In this paper we shall consider countably infinite locally $2$-arc-transitive bipartite graphs, which, by the comments above, include as a subclass the 
countably infinite locally $2$-arc-transitive graphs for which $G$ is intransitive on vertices. Our starting point is work of Droste \cite{Droste1} on 
\emph{semilinear orders} and later of Warren (and others) on the, so called, $3$-$CS$-transitive \emph{cycle-free partial orders}. A surprising 
by-product of that work was that it led to a new and interesting family of countably infinite locally $2$-arc-transitive bipartite graphs, essentially 
constructed by gluing linear orders together in a prescribed way (see Section~2 for a detailed description of this construction). This family of 
examples is of particular significance since it has no obvious finite analogue, that is, they are not simply examples that arise by taking some known 
finite family and allowing a parameter to go to infinity. Motivated by this, recent attempts have been made to see to what extent the methods of 
Warren may be used to start to try and understand countably infinite locally $2$-arc-transitive bipartite graphs. Some recent work relating to these 
attempts includes \cite{Gray1} and \cite{Gray2}. However, until now the results have been fairly negative, highlighting some of the difficulties 
involved in generalizing Warren's techniques to obtain bipartite graphs with the required level of transitivity. In this article we succeed in an 
extension of the methods that does give (many) new examples of countably infinite locally $2$-arc-transitive bipartite graphs, and in the process we 
discover some interesting connections with semilinear spaces (satisfying certain point--line--point transitivity properties) and the highly 
arc-transitive digraphs of Cameron et al. \cite{Cameron2}. 

The paper is structured in the following way. In Section~2 we give the basic definitions and ideas required for the rest of the paper, and we explain 
our approach to the study of infinite locally $2$-arc-transitive bipartite graphs. The idea is to reduce the analysis of such bipartite graphs to that 
of families determined by a poset derived from the bipartite graph, which we call its \emph{interval}. We go on in Section~3 to look at one of these 
families in detail, namely the one in which the intervals are totally ordered, and the main results and constructions of the paper are given in that 
section. We concentrate on presenting a construction initially based on the integers, which is generalized to many of the countable 1-transitive 
linear orders in Morel's list \cite{Morel}, as well as one `2-coloured' case. As we shall explain, the general case gives rise to many more 
complications, even under the restrictions we have chosen to adopt. The special case treated is that in which the interval has many pairs of 
consecutive points, which can themselves be viewed as edges of a digraph $\Delta$, whose role is in some sense to measure how far our original graph 
is from being cycle-free. So a dense linear order $\mathbb Q$, or its 2-coloured analogue ${\mathbb Q}_2$ clearly cannot arise in this setting, and to 
cover this situation, or indeed any in which the two levels of $\Delta$ are embedded in non-consecutive levels seems technically considerably more 
complicated. We therefore content ourselves with making some remarks about the extension of these methods to a more general situation at the end of 
section 3 without full details. 

The main results of the paper are Theorems \ref{3.15} and \ref{3.16}. The link between bipartite graphs and posets, which lies at the heart of the 
paper, is described in Theorem \ref{3.15}, which gives the correspondence between connected highly arc-transitive digraphs fulfilling some natural 
conditions and countable connected sufficiently transitive bipartite graphs. Theorem \ref{3.16} builds on this, and asserts the existence of a 2-level 
partial order $P$ corresponding to certain possible `inputs', consisting of a choice of $\Delta$ and a countable 1-transitive linear order $\cal Z$. 
These 2-level partial orders (meaning that all maximal chains have length 2) can alternatively be construed as bipartite graphs, whose parts are the 
sets of minimal and maximal points respectively, and the ones constructed are then locally 2-arc transitive. The theorem characterizes which such 
partial orders can arise from the construction, in terms of the values of $\Delta$ and $\cal Z$, demonstrating that this constitutes a substantial 
generalization of the cycle-free situation. 

Some results concerning arbitrary non-totally ordered intervals are given in Section~4, where we treat some special cases. Since these do not fit in 
with the other main themes of the paper, and the methods are rather ad hoc, we give rather few details.

\section{Basic definitions and outline of approach}

\subsection{Definitions and background}

For clarity and completeness we first recall the precise definitions of `graph', `directed graph', and `partially ordered set'. A {\em graph} $\Gamma$
is formally a pair of sets $(V\Gamma, E\Gamma)$ where $E\Gamma$ is a set of 2-element subsets of $V\Gamma$ (thus we restrict attention to `simple' 
graphs, those without loops or multiple edges), and a {\em directed graph} or {\em digraph}, is a pair $(V\Gamma, E\Gamma)$ of sets for which 
$E\Gamma$ is a set of ordered pairs of distinct members of $V\Gamma$, and such that for no $x$ and $y$ are both $(x, y)$ and $(y, x)$ edges. (This is 
called `asymmetric', also `oriented graph'; it is sometimes, for instance in \cite{Cherlin}, required in the definition of `digraph', but not 
always. Here it simplifies the exposition since all our digraphs are asymmetric anyhow.) We also use the notation $(D, \to)$ for a digraph, so that 
$x \to y$ is an alternative way of saying that $(x, y)$ is a directed edge. A {\em partially ordered set} (or {\em poset}) is a pair $(M, <)$ such 
that $<$ is a binary relation on $M$ which is irreflexive ($\neg x < x$ for all $x \in M$) and transitive. We may also work with the corresponding 
reflexive partial order $\le$ on $M$ defined by $x \le y \Leftrightarrow x < y$ or $x = y$. 

We begin proceedings with an illustrative example, to give the reader a feel for the ideas that will be formalized below.

The key idea exploited in much of the paper is the interplay between partial orders and directed graphs (often bipartite). Thus it is an apparent 
triviality that a bipartite graph with parts $X$ and $Y$ may be viewed as a partial order by selecting $X$ say as the `lower level' and writing 
$x < y$ for $x \in X$ and $y \in Y$ if $x$ and $y$ are adjacent. It is trivial in the sense that transitivity is vacuously true. One can however form 
an extension of any partial order $(M, <)$, termed its `Dedekind--MacNeille completion', written $M^D$ (or $(M^D, <)$) which is a partial order 
analogue of the familiar Dedekind-completion of a linear order. This may have additional points, and the 2-level partial order just mentioned may have 
many extra points `in between', whose structure will throw light on the disposition of the points in the original graph or digraph.

To explain our first example, we need the following basic definitions, which may also be found in \cite{Warren,Droste4}.

A partial order $(M, <)$ is said to be {\em Dedekind--MacNeille complete} ({\em $DM$-complete} for short) if any maximal chain is Dedekind-complete, and any 
two-element subset bounded above (below) has a least upper bound or supremum (a greatest lower bound or infimum respectively). This is equivalent to 
saying that any non-empty bounded above (below) subset has a supremum (infimum respectively). It is shown for instance in \cite{Davey} that for any 
partially ordered set $(M, <)$ there is a minimal $DM$-complete partial order $(M^D, <)$ containing it, which is unique up to isomorphism, and this is 
called its {\em Dedekind-MacNeille completion}. An efficient way to describe $M^D$ explicitly, which we shall sometimes need, is as the family of 
`ideals' of $M$, where an {\em ideal} is a non-empty bounded subset $J$ of $M$ which is equal to the set of lower bounds of its set of upper bounds, 
partially ordered by inclusion. The original set $M$ is then identified with the family of principal ideals, being subsets of $M$ of the form 
$\{x \in M: x \le a\} = M^{\le a}$ for some $a \in M$, and $M$ is already $DM$-complete provided that all ideals are principal. Note that the notions 
of Dedekind-MacNeille completion in \cite{Davey} and \cite{Warren} are slightly different---in the former case it is a complete lattice, so contains a 
least and greatest (obtained by deleting the restrictions that an ideal be non-empty and bounded), but not in the latter. Since the two notions differ 
at most on these two points, we regard this as a minor problem.

To illustrate the definition in simple cases, first note that the Dedekind--MacNeille completion of the rationals $\mathbb Q$ is the reals $\mathbb R$,
since the ideals are precisely left Dedekind cuts of the form $(-\infty, a] = \{x \in {\mathbb Q}: x \le a\}$ for some real number $a$ (its set 
of upper bounds equals $[a, \infty)$, whose set of lower bounds is $(-\infty, a]$ again. The version as in \cite{Davey} would be 
${\mathbb R} \cup \{\pm \infty\}$). For a finite example, we give in Figure 1(b) the Dedekind--MacNeille completion of the partial order of Figure 
1(a). The two non-principal ideals are $p = \{x, y\}$ and $q = \{y, z\}$. For instance, the set of upper bounds of $\{x, y\}$ equals
$\{u, v\}$, whose set of lower bounds is just $\{x, y\}$ again.

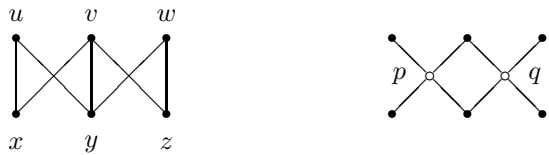
\begin{figure}[h]
\unitlength=1mm
\begin{picture}(120.00,20.00)

\put(40,5){\circle*{1.2}}
\put(50,5){\circle*{1.2}}
\put(60,5){\circle*{1.2}}
\put(40,15){\circle*{1.2}}
\put(50,15){\circle*{1.2}}
\put(60,15){\circle*{1.2}}

\put(40,5){\line(0,1){10}}
\put(50,5){\line(0,1){10}}
\put(60,5){\line(0,1){10}}
\put(40,5){\line(1,1){10}}
\put(50,5){\line(-1,1){10}}
\put(50,5){\line(1,1){10}}
\put(60,5){\line(-1,1){10}}

\put(40,1){\makebox(0,0)[cc]{$x$}}
\put(50,1){\makebox(0,0)[cc]{$y$}}
\put(60,1){\makebox(0,0)[cc]{$z$}}
\put(40,18){\makebox(0,0)[cc]{$u$}}
\put(50,18){\makebox(0,0)[cc]{$v$}}
\put(60,18){\makebox(0,0)[cc]{$w$}}

\put(90,5){\circle*{1.2}}
\put(100,5){\circle*{1.2}}
\put(110,5){\circle*{1.2}}
\put(95,10){\circle{1.2}}
\put(105,10){\circle{1.2}}
\put(90,15){\circle*{1.2}}
\put(100,15){\circle*{1.2}}
\put(110,15){\circle*{1.2}}

\put(90,5){\line(1,1){4.5}}
\put(100,5){\line(-1,1){4.5}}
\put(100,5){\line(1,1){4.5}}
\put(110,5){\line(-1,1){4.5}}
\put(90,15){\line(1,-1){4.5}}
\put(100,15){\line(-1,-1){4.5}}
\put(100,15){\line(1,-1){4.5}}
\put(110,15){\line(-1,-1){4.5}}

\put(91,10){\makebox(0,0)[cc]{$p$}}
\put(109,10){\makebox(0,0)[cc]{$q$}}

\end{picture}

\caption{(a) a finite partial order and (b) its Dedekind--MacNeille completion}
\end{figure}

We say that $(M, <)$ is {\em cycle-free} if in $M^D$ there is a unique path between any two points of $M$. Note that any cycle-free partial order is 
automatically `connected', which means that there is {\em at least one} path (in $M^D$) between any two points. We usually write `cycle-free partial 
order' as $CFPO$. Here, `path' is thought of in a graph-theoretical sense, being a finite union of linear segments $[a, b]$ where $a < b$ which are 
maximal chains in the induced partial order of points between $a$ and $b$, and such that the segments only overlap at endpoints, and its {\em length} 
is the number of such segments. Note that in Figure 1(a), $x, u, y, v, x$ `looks like' a cycle, but should not count as such, since in the completion, 
the corresponding path has duplicated sections between $x$ and $p$ for instance (this decision concerning the correct definition of `cycle-freeness' 
is forced on us by the evident requirement that any substructure of a cycle-free partial order should be cycle-free). Despite this, the partial order 
does {\em not} count as cycle-free, but this is because there is more than one path {\em in the completion} between $y$ and $v$. 

A full development of the theory of cycle free partial orders 
is given in \cite{Warren} but we remark on the most important features. First, it follows from cycle-freeness that the interval 
$[a, b] = \{x \in M^D: a \le x \le b\}$ {\em is} a linear order (even though in the definition, we had to say we were taking a maximal chain in this 
interval). Next, we have clearly defined notions of `cone' and `ramification order' which are given as follows. For any partially ordered set $M$, we 
say that $x \in M^D$ is an {\em upward ramification point of $M$} ({\em downward ramification point of $M$}) if there are incomparable points 
$a, b \in M$ such that $x$ is the infimum (supremum respectively) of $a$ and $b$, written $x = a \wedge b$ ($x = a \vee b$ respectively). Intuitively, 
these are points at which the structure branches going up and down respectively. 
For instance, if $M$ is the poset given in Figure 1(a) then the points $p, q \in M^D$ in Figure 1(b) are each both upward and downward ramification points. 
We write $\uparrow \hspace{-.05in} \mathrm{Ram}(M)$ and 
$\downarrow \hspace{-.05in} \mathrm{Ram}(M)$ for the sets of upward and downward ramification points, respectively, of $M$, $\mathrm{Ram}(M)$ for 
$\uparrow \hspace{-.05in} \mathrm{Ram}(M) \hspace{.05in} \cup \downarrow \hspace{-.05in} \mathrm{Ram}(M)$ and $M^+$ for the union of $M$ and 
$\mathrm{Ram}(M)$. Note that $\uparrow \hspace{-.05in} \mathrm{Ram}(M)$ and $\downarrow \hspace{-.05in} \mathrm{Ram}(M)$ are each subsets of $M^D$. In 
general $\uparrow \hspace{-.05in} \mathrm{Ram}(M) \neq \hspace{.05in} \uparrow \hspace{-.05in} \mathrm{Ram}(M^D)$, but for cycle-free partial orders 
the two sets coincide; see \cite{Warren} Lemma~2.3.11 and also Lemma \ref{3.7} below. An {\em upper cone} of a point $x$ is an equivalence class under 
the relation $\sim$ given by $a \sim b$ if there is $y > x$ such that $a, b \ge y$ (which is an equivalence relation under the assumption of 
cycle-freeness, or even under the weaker assumption, explored later in the paper, that all intervals of $M^D$ are chains), and the number of upper 
cones is called the (upper) {\em ramification order} of $x$ (similarly for lower ramification order). If the partial ordering is finite, or even 
discrete, each cone will be represented by a minimum element, and then we can view the ramification order as the number of branches at that point in a 
direct sense; the version involving cones is a way of expressing this when discreteness may fail. Then (still in the cycle-free case) $M^+$ is the 
smallest extension of $M$ in which all 2-element subsets of $M$ have least upper bounds and greatest lower bounds (provided they are correspondingly 
bounded in $M$), and a useful property of $M^+$ (over $M^D$) is that (in the infinite case) its cardinality is no greater than that of $M$. We shall 
see in Lemma \ref{3.7} that the same applies in any `diamond-free' partial order, so at that point we shall be free to work with $M^+$ rather than 
$M^D$. If $a$ is an upward ramification point then its {\em ramification order}, written $\uparrow \hspace{-.05in} ro(a)$ is the number of upper cones 
at $a$ (similarly $\downarrow \hspace{-.05in} ro(a)$ at a downward ramification point). If these are the same for all ramification points then we may 
write $\uparrow \hspace{-.05in} ro(M)$ or $\downarrow \hspace{-.05in} ro(M)$.

\subsection{A cycle-free example and a lemma}
\label{subsec_CFPO}

The simplest example of a cycle-free partial order from \cite{Warren} which gives rise to a locally 2-arc-transitive bipartite graph is as follows, 
and this illustrates the general method. We start with the `inside' of the structure, which can be viewed as the `scaffolding', to be thrown away at 
the end, and this is taken to be an $(m,n)$-directed tree $T$ for any $m$ and $n$ which are at least 2 and at most $\aleph_0$. This is a partial order 
which arises as the transitive closure of a certain digraph, such that all vertices have in-valency $m$ and out-valency $n$, and subject to this is 
`free'. An explicit construction of such a digraph is done in countably many stages. The first stage consists of a single copy of $\mathbb Z$ under 
its usual (linear) ordering. At the next, for each point of $\mathbb Z$ we add new points so that each original point has $m-1$ new points immediately 
above it and $n-1$ new points immediately below it, and each upper point has an infinite ascending sequence above it and each lower point has an 
infinite descending sequence below it. Thus all the original points have the correct valencies, but the new ones do not. At a general stage, all the 
points from the previous stage which do not yet have the correct valencies are assigned new points above or below, and sequences above them and below 
them. Then we take the union, and the result is a (directed) tree, since none of the new points are identified, so we never get any circuits.

The next step is to choose a countable dense set of maximal chains of $T$, and put a point above and below each such chain. This can be done 
explicitly by choosing for each point of $T$ a maximal chain passing through that point. The result is partially ordered by saying that $t \in T$ lies 
below the point $x$ above a maximal chain $C$ provided $t$ is below some point of $C$ (and similarly for $t$ lying above the point below $C$). We let 
$M$ be the family of all these choices, and this receives the induced partial ordering. Since we have now `thrown away' all the points of $T$, it is 
clear that $M$ is a 2-level partial order (meaning that all its maximal chains have length 2), so may be regarded as a bipartite graph (or digraph). 
The points of $T$ are not however completely lost, as one checks that $M \cup T$ is equal to the Dedekind--MacNeille completion $M^D$ of $M$. In 
Figure 2 we illustrate the second stage in the construction of $T$ for the case $m = 3$, $n = 2$ (where for convenience we put the in-neighbours to 
the right and the out-neighbours to the left, though this has no significance).

\begin{figure}[h]
\unitlength=1mm
\begin{picture}(120.00,100.00)

\put(75,20){\circle*{1.2}}
\put(75,30){\circle*{1.2}}
\put(75,40){\circle*{1.2}}
\put(75,50){\circle*{1.2}}
\put(75,60){\circle*{1.2}}
\put(85,60){\circle*{1.2}}
\put(95,60){\circle*{1.2}}
\put(85,70){\circle*{1.2}}
\put(95,70){\circle*{1.2}}
\put(85,60){\circle*{1.2}}
\put(95,60){\circle*{1.2}}
\put(85,50){\circle*{1.2}}
\put(95,50){\circle*{1.2}}
\put(85,40){\circle*{1.2}}
\put(95,40){\circle*{1.2}}
\put(85,30){\circle*{1.2}}
\put(95,30){\circle*{1.2}}
\put(65,50){\circle*{1.2}}
\put(65,40){\circle*{1.2}}
\put(65,30){\circle*{1.2}}
\put(65,20){\circle*{1.2}}
\put(65,10){\circle*{1.2}}

\put(75,20){\line(0,1){10}}
\put(75,30){\line(0,1){10}}
\put(75,40){\line(0,1){10}}
\put(75,50){\line(0,1){10}}
\put(75,17){\line(0,1){3}}
\put(75,60){\line(0,1){3}}
\put(75,50){\line(1,1){10}}
\put(65,50){\line(1,1){10}}
\put(65,40){\line(1,1){10}}
\put(65,30){\line(1,1){10}}
\put(65,20){\line(1,1){10}}
\put(65,10){\line(1,1){10}}
\put(75,40){\line(1,1){10}}
\put(75,30){\line(1,1){10}}
\put(75,20){\line(1,1){10}}
\put(75,40){\line(2,1){20}}
\put(75,30){\line(2,1){20}}
\put(75,20){\line(2,1){20}}
\put(75,50){\line(2,1){20}}
\put(85,60){\line(0,1){10}}
\put(95,60){\line(0,1){10}}
\put(85,70){\line(0,1){3}}
\put(95,70){\line(0,1){3}}
\put(85,50){\line(0,1){2}}
\put(95,50){\line(0,1){2}}
\put(85,40){\line(0,1){2}}
\put(95,40){\line(0,1){2}}
\put(85,30){\line(0,1){2}}
\put(95,30){\line(0,1){2}}
\put(63,46){\line(1,2){2}}
\put(63,36){\line(1,2){2}}
\put(63,26){\line(1,2){2}}
\put(63,16){\line(1,2){2}}
\put(63,6){\line(1,2){2}}

\put(75,69){\line(0,1){1}}
\put(75,72){\line(0,1){1}}
\put(85,79){\line(0,1){1}}
\put(85,82){\line(0,1){1}}
\put(95,79){\line(0,1){1}}
\put(95,82){\line(0,1){1}}
\put(75,10){\line(0,1){1}}
\put(75,7){\line(0,1){1}}

\end{picture}

\caption{Second stage in the construction of $T$}
\end{figure}
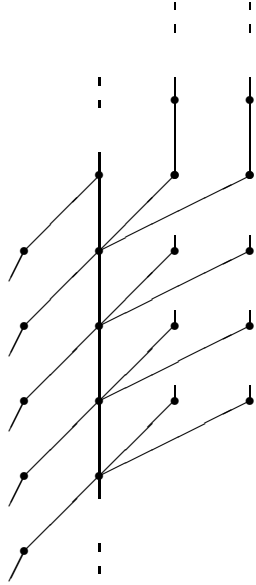

There is a straightforward link between a natural notion of transitivity on a 2-level $CFPO$, and its arc-transitivity when viewed as a directed 
graph. For $k \ge 1$ the notions introduced in \cite{Warren} are as follows. We say that $M$ is {\em $k$-connected-set-homogeneous} 
($k$-$CS$-homogeneous for short) if any isomorphism between {\em connected} $k$-element substructures extends to an automorphism. (Without the 
`connectedness' stipulation, the notion is essentially vacuous, since for instance there are 2-element antichains which are at different `distances' 
in the partial order, and no automorphism could take one to the other.) 
This is the model theorist's notion of `substructure', corresponding to the graph theorist's notion of `induced substructure'.
It is {\em $k$-$CS$-transitive} if for any two isomorphic connected $k$-element 
substructures there is an automorphism taking one to the other (not necessarily extending an originally given isomorphism). 
A bipartite digraph is 
{\em $2$-arc-transitive} if its underlying graph is 2-arc-transitive, and similarly for local 2-arc-transitivity. (Note the non-standard convention we 
adopt here; usually a digraph would be called $2$-arc-transitive if the automorphism group acted transitively on directed $2$-arcs.) It is then clear 
that local 2-arc-transitivity is equivalent to 3-$CS$-homogeneity. The 2-arcs are of two types, shaped like $\Lambda$s or Vs, so 3-$CS$-transitivity 
just says that the automorphism group acts transitively on each of the families of $\Lambda$s and of Vs. We follow the convention that `local 
2-arc-transitivity' is used when talking of a graph or digraph, and `3-$CS$-homogeneity' is used when talking of a partial order. 

The proof of the following lemma is straightforward and is omitted. 

\begin{lemma}\label{2.1}
Let $M$ be a locally $1$- and $2$-arc-transitive bipartite graph, viewed as a $2$-level partial order $M = X \cup Y$ with $X$ and $Y$ the sets of 
minimal and maximal points respectively. Let $G = {\rm Aut}(M)$. Then

(i) for all $x, x' \in X$, $y,y' \in Y$ such that $x < y$ and $x' < y'$, we have $[x,y] \cong [x',y']$,

(ii) $G$ acts transitively on $\downarrow \hspace{-.05in} \mathrm{Ram}(M)$, and transitively on $\uparrow \hspace{-.05in} \mathrm{Ram}(M)$,

(iii) $\uparrow \hspace{-.05in} \mathrm{ro}(a) = \hspace{.05in} \uparrow \hspace{-.05in} \mathrm{ro}(b)$ for all $a,b \in \hspace{.05in} \uparrow 
\hspace{-.05in} \mathrm{Ram}(M)$ and $\downarrow \hspace{-.05in} \mathrm{ro}(a) = \hspace{.05in} \downarrow \hspace{-.05in} \mathrm{ro}(b)$ for all 
$a,b \in \hspace{.05in} \downarrow \hspace{-.05in} \mathrm{Ram}(M)$, in cases where these are defined.
\end{lemma}

In part (i) of this lemma, $[x, y]$ `officially' stands for $[x, y]^{M^D} = \{z \in M^D: x \le z \le y\}$. In practice however it is much easier to 
work with $[x, y]^{M^+}$, since this is countable provided $M$ is, and the lemma is true for this set as well. We call the poset $[x, y]^{M^D}$ the 
\emph{interval} of the bipartite graph $M$, and we denote it by $I(M)$, and $[x, y]^{M^+}$ is written $I^+(M)$. Strictly speaking, we should really 
regard $I^+(M)$ as a `2-coloured' chain, meaning that every point apart from the endpoints has at least one of 2 colours assigned, corresponding to 
whether it is an upward or downward ramification point. Once again, part (i) of the lemma holds even if the isomorphism is required to preserve 
colours. In the simplest case, all points apart from the endpoints are coloured by both colours. Viewing $M$ as a partial order from the beginning, 2- 
and 3-$CS$-transitivity are sufficient hypotheses for the lemma to apply. We remark that in Warren's work the focus is slightly wider, namely on 
3-$CS$-transitivity rather than 3-$CS$-homogeneity. Since we principally wish to explore the connection with 2-arc-transitivity, we restrict to the 
stronger hypothesis of 3-$CS$-homogeneity (which is already extremely rich).

\subsection{Outline of approach}\label{sec_outline}

The results of Warren and others on cycle-free partial orders suggest the following approach for the analysis of countable locally $2$-arc-transitive 
bipartite graphs. 

(I) Identify the possible partial orders that can arise as intervals $I(M)$ where $M$ is a connected countable locally 2-arc-transitive bipartite 
graph. 

(II) For each possible interval $I$ (or class of intervals) classify those connected countable locally 2-arc-transitive bipartite graphs such that 
$I(M) \cong I$.  

Of course, we do not mean to suggest that we believe there is any chance that a complete classification of countable 2-arc-transitive bipartite graphs 
is really achievable, but this approach does at least break the problem down in such a way that for certain intervals, or families of intervals, one 
may obtain quite a lot of information about the corresponding class. It is (II) which presents the main problem, as even for any of the possible 
infinite (2-coloured) linear orders which can arise as $I(M)$, there will be many many possibilities for $M$, as we shall remark in section 3. The 
answer to (I) in the case of linear intervals is given in Lemma \ref{3.9}. 

From the definitions, if a poset arises as an interval of a locally 2-arc-transitive bipartite graph then it must itself be $DM$-complete, and it 
will have to have a maximal and minimal element, but as we shall see below, not every $DM$-complete poset with a maximal and minimal element will 
arise in this way. On the other hand, the construction given in \cite{Gray2} demonstrates that if the condition on $M$ is weakened to local 
1-arc-transitivity then any $DM$-complete poset with a maximal and minimal element can arise.

\section{Completions with chain intervals}

The motivating examples for the approach outlined in the previous section are the 3-$CS$-homogeneous 2-level partial $CFPO$s arising from Warren's 
classification. All of these examples have the property that $I(M)$ is a chain. Hence a natural next step is to consider those locally 
2-arc-transitive bipartite graphs $M$ such that $I(M)$ is a chain. 

This is equivalent to saying that the completion $M^D$ does not embed any \emph{diamonds}, where a diamond is a poset with elements $\{ a,b,c,d \}$ 
such that the only non-trivial relations are $a < b < d$ and $a < c < d$. So it would not be unreasonable to refer to this as the class of 
\emph{diamond-free partial orders}. The analysis of this family of posets splits naturally into cases depending on whether or not the interval is 
finite. First we show that it makes no difference to require the intervals of $M^D$ or of $M^+$ to be chains.

\begin{lemma}\label{3.1} If $M$ is a partial order in which for every $a \le b$, $[a, b]^{M^+}$ is a chain, then for any $a \le b$ in $M^D$, 
$[a, b]^{M^D}$ is a chain.    \end{lemma}

\noindent{\bf Proof}: As remarked in section 2, the standard construction of $M^D$ is as the family of ideals of $M$ ordered by inclusion, where an 
`ideal' is a non-empty bounded above subset $J$ which is equal to the set of lower bounds of its set of upper bounds. Thus if we write $J^\uparrow$ 
and $J^\downarrow$ for the set of upper, lower bounds of $J$, respectively, then the condition says that $J, J^\uparrow \neq \emptyset$, and 
$J^{\uparrow \downarrow} = J$. Also $M$ is embedded in $M^D$ by the map sending $a$ to $M^{\le a} = \{x \in M: x \le a\}$. Thus if $J \in M^D$, and $a 
\in M$, $a \le J \Leftrightarrow a \in J$ under this identification.

Now suppose for a contradiction that $x, y$ are incomparable members of $[a, b]$ in $M^D$. Thus $a$ and $b$ are ideals, hence bounded above and below 
in $M$, so by decreasing $a$ and increasing $b$ if necessary, we may suppose that $a, b \in M$. Viewing $x$ and $y$ as subsets of $M$ (ideals), since 
they are incomparable, $x \not \subseteq y$ and $y \not \subseteq x$, so we may choose $u \in x \setminus y$ and $v \in y \setminus x$. Thus $u, v \in 
M$, and using the above convention, $u \le x$, $u \not \le y$, $v \le y$, and $v \not \le x$. Since $a, u \le x$, $a \vee u$ exists and lies in $M^+$ 
and similarly, $a \vee v \in M^+$, and clearly, $a \le a \vee u, a \vee v \le b$. Since we are assuming that $[a, b]^{M^+}$ is a chain, $a \vee u \le a 
\vee v$ or $a \vee v \le a \vee u$, assume the former. But then $u \le a \vee u \le a \vee v \le y$, which is a contradiction.   $\Box$

\vspace{.1in}

We need the following result from \cite{Warren}, (the `density lemma', 2.4.7) on several occasions. There it was proved just for $CFPO$s, but here we 
need it in the more general situation that all intervals are linear. As in the above proof, for a subset $a$ of $M$ we write $a^\uparrow$ for 
$\{x \in M: a \le x\}$ and $a^\downarrow$ for $\{x \in M: a \ge x\}$.

\begin{lemma}\label{3.2} For any partial order $M$ in which for every $x < y$, $[x, y]^{M^+}$ is a chain, if $a < b$ in $M^D$, then both $[a,b) \cap 
(M \cup \uparrow \hspace{-.05in} {\rm Ram}(M))$ and $(a,b] \cap (M \cup \downarrow \hspace{-.05in} {\rm Ram}(M))$ are non-empty.
\end{lemma}

\noindent{\bf Proof}: By Lemma \ref{3.1}, all intervals of $M^D$ are also chains. 

Viewing $a$ and $b$ as ideals, $a \subset b$. Since any ideal is non-empty and bounded, there are $x, y \in M$ such that $x \le a < b \le y$. Clearly 
$b^\uparrow \subseteq a^\uparrow$. By definition of `ideal', $a^{\uparrow \downarrow} = a$ and $b^{\uparrow \downarrow} = b$, so $b^\uparrow = 
a^\uparrow$ would imply that $a = a^{\uparrow \downarrow} = b^{\uparrow \downarrow} = b$, a contradiction. Hence $b^\uparrow \subset a^\uparrow$, and 
there is $z \in a^\uparrow \setminus b^\uparrow$. Hence $a \le z$ and $b \not \le z$. In $M^D$ this becomes $a \subseteq M^{\le z}$ and $b \not \subseteq 
M^{\le z}$. If $M^{\le z} \subset b$ then $M^{\le z} \in [a, b) \cap M$. Otherwise, $y$ and $z$ are incomparable, as $y \not \le z$ since $b \not \le 
z$ and if $z \le y$ then $b, z \in [x, y]$, which is assumed linear, which implies $z \le b$ after all. Hence $y \wedge z$ exists and lies in 
$\downarrow \hspace{-.05in} \mathrm{Ram}(M)$. Also, $y \wedge z \in [a, b)$ since it lies in $[x, y]$ which is linear, so is comparable with $b$, and 
the only option is $y \wedge z < b$. 

For the other clause, as $a \subset b$, there is $t \in b \setminus a$. As $b = b^{\uparrow \downarrow}$, $t \le b^\uparrow$, and similarly, $t \not \le 
a^\uparrow$. Therefore $t \le b$ and $t \not \le a$. If $t = b$ then $t \in (a, b] \cap M$. Otherwise, $t < b$ and by a similar argument to the first 
clause, $x \vee t \in [a, b) \cap \uparrow \hspace{-.05in} \mathrm{Ram}(M)$. $\Box$

\vspace{.1in}

\begin{assumption}\label{3.3}
Throughout the rest of this section $M$ will be a 2-level connected, countable partial order such that $I(M)$ is a chain. 
\end{assumption}

\subsection{Finite chain diamond-free partial orders}
\label{subsec_finitechain}

When the interval $I(M)$ is a finite chain, there is a natural connection with incidence structures, specifically with linear, and semilinear spaces (which are also called partial linear spaces in the literature).   

\begin{definition}\label{3.4}
A \emph{semilinear space} is a pair $(P, L)$ consisting of a set $P$ of \emph{points} and a set $L$ of \emph{lines} such that

(i) every line contains at least two points, and

(ii) any two distinct points are on at most one line.
\end{definition}

The following proposition gives an immediate connection of this notion with the types of partial order we are considering.

\begin{prop}\label{3.5} Let $M$ be a bipartite graph viewed as a $2$-level partial order $M = X \cup Y$ where $X$ and $Y$ are the sets of minimal and 
maximal points respectively, both non-empty, and such that any member of $Y$ is above at least two points of $X$. Then $M$ is $DM$-complete if and 
only if $M$ is the incidence graph of a semilinear 
space. 
\end{prop}
\noindent{\bf Proof}: First suppose that $M$ is the incidence graph of a semilinear space. The maximal chains of $M$ are clearly all 
Dedekind-complete. Therefore, to show that $M$ is $DM$-complete it suffices to show that any two-element subset $A$ of $M$ bounded above (below) has a 
supremum (infimum). If $A \subseteq X$, this follows from the second clause in Definition \ref{3.4}. If $A \subseteq Y$ is bounded below by $x_1$ and 
$x_2$ in $X$, then $\{x_1, x_2\}$ has more than one upper bound in $Y$, again contrary to the same clause. Conversely, assume that $M$ is 
$DM$-complete. Clause (i) of Definition \ref{3.4} follows from the stated hypothesis. To verify (ii), let $x_1, x_2 \in X$ be distinct points lying on 
a line. Then $\{x_1, x_2\}$ is bounded above, so its supremum exists, and this means that the line that $x_1$ and $x_2$ lie on is unique. $\Box$

\vspace{.1in}

Now introducing transitivity hypotheses we have the following.

\begin{prop}\label{3.6}
Let $M$ be a connected locally $2$-arc-transitive bipartite graph, viewed as a $2$-level partial order $M = X \cup Y$ with $X$ and $Y$ the sets of 
minimal and maximal points respectively. If $I = I(M)$ is a finite chain then $|I| = 2$ or $|I| = 3$. Moreover:

(i) If $|I| = 3$ then $M$ is complete bipartite.

(ii) If $|I| = 2$ and $|X|, |Y| \ge 2$ then $M$ is the incidence graph of a semilinear space $S$ such that ${\rm Aut}(S)$ is transitive on 
configurations of the form $(p_1,l,p_2)$, where $l$ is a line incident with $p_1$  $p_2$, and on configurations $(l_1,p,l_2)$ where $p$ is a point 
incident with $l_1$ and $l_2$.     
\end{prop}
\noindent{\bf Proof}:
Let $a \in X$, $b \in Y$ with $a < b$, and suppose for a contradiction that $|[a, b]| \ge 4$. Let $a < x < y < b$. Since $a < x$ there is $c \neq a$ 
in $X$ with $c < x$ and since $x < y$ there is $d \in X$ with $d < y$ and $d \not < x$. Then $a \vee c, a \vee d \in [a, b]$ so as this is a chain, $a 
\vee c < a \vee d$ or $a \vee c \ge a \vee d$. But the latter would imply that $d < x$, contrary to its choice. Hence $a \vee c < a \vee d$. By 3-$CS$
-transitivity there is an automorphism taking $\{a, c, b\}$ to $\{a, d, b\}$. This must take $a \vee c$ to $a \vee d$ and fix $b$, so its orbit is an 
infinite subset of $[a, b]$, contrary to $I$ finite. See Figure 3(a).

For (i), let $a < b$ and $[a, b] = \{a, x, b\}$. Then by Lemma \ref{3.2} $x$ must be an upward and downward ramification point, so by using 
3-$CS$-transitivity we see that all ramification points are on the middle level. We show that every element of $X$ lies below $x$. Suppose not, for a 
contradiction. Choose $c \in X$ not below $x$. As $M$ is connected, there is a path from $c$ to $a$. Choose $c$ so that this path is of least possible 
length, and let it begin $c < b_1 > a_1$. By minimality of the length of the path, $a_1 < x$. Then $b, b_1 > a_1$, so $b \wedge b_1$ exists and is the 
middle element of $[a_1, b] = \{a_1, x, b\}$. Therefore $b \wedge b_1 = x$, and so $x < b_1$. From $a_1, c < b_1$ it follows that $a_1 \vee c$ exists 
and is the middle element of $[a_1, b_1] = \{a_1, x, b_1\}$. Hence $a_1 \vee c = x$, and so $c < x$ after all, giving a contradiction. See Figure 3(b).

Similarly, every element of $Y$ lies above $x$, and it follows that every element of $X$ is below every element of $Y$, so that $M$ is complete 
bipartite.

For (ii), we take the members of $X$ as `points' and those of $Y$ as `lines'. Saying that every line contains at least two points then just says that 
every member of $Y$ is above at least 2 members of $X$. Now as $M$ is connected and $X$ and $Y$ each have size at least 2, there are upward and 
downward ramification points, and as $|I| = 2$, the upward ramification points must be the points of $X$ and the downward ramification points are the 
points of $Y$. (The fact that {\em all} members of $X$, $Y$ are ramification points follows from 3-$CS$-transitivity.) It follows that every member of 
$Y$ has at least 2 members of $X$ below it. Since $|I| = 2$, $M$ is $DM$-complete, so by Proposition \ref{3.5}, $M$ is the incidence graph of a 
semilinear space $S$. The transitivity properties of ${\rm Aut}(S)$ stated amount to a reformulation of 3-$CS$-homogeneity.  $\Box$

\begin{figure}
\unitlength=1mm
\begin{picture}(120.00,20.00)

\put(35,2){\circle*{1.2}}
\put(40,2){\circle*{1.2}}
\put(45,2){\circle*{1.2}}
\put(35,7){\circle{1.2}}
\put(35,12){\circle{1.2}}
\put(35,17){\circle*{1.2}}

\put(35,2){\line(0,1){4.5}}
\put(40,2){\line(-1,1){4.5}}
\put(45,2){\line(-1,1){9.5}}
\put(35,7.5){\line(0,1){4}}
\put(35,12.5){\line(0,1){4.5}}

\put(21,10){\makebox(0,0)[cc]{(a)}}
\put(32,2){\makebox(0,0)[cc]{$a$}}
\put(38,1){\makebox(0,0)[cc]{$c$}}
\put(48,2){\makebox(0,0)[cc]{$d$}}
\put(32,7){\makebox(0,0)[cc]{$x$}}
\put(32,12){\makebox(0,0)[cc]{$y$}}
\put(32,18){\makebox(0,0)[cc]{$b$}}

\put(66,10){\makebox(0,0)[cc]{(b)}}

\put(80,5){\circle*{1.2}}
\put(85,5){\circle*{1.2}}
\put(90,5){\circle*{1.2}}
\put(80,10){\circle{1.2}}
\put(80,15){\circle*{1.2}}
\put(90,15){\circle*{1.2}}

\put(80,5){\line(0,1){4.5}}
\put(85,5){\line(-1,1){4.5}}
\put(80,10.5){\line(0,1){4.5}}
\put(85,5){\line(1,2){5}}
\put(90,5){\line(0,1){10}}

\put(80,2){\makebox(0,0)[cc]{$a$}}
\put(85,2){\makebox(0,0)[cc]{$a_1$}}
\put(90,2){\makebox(0,0)[cc]{$c$}}
\put(77,10){\makebox(0,0)[cc]{$x$}}
\put(80,18){\makebox(0,0)[cc]{$b$}}
\put(90,18){\makebox(0,0)[cc]{$b_1$}}

\put(105,5){\circle*{1.2}}
\put(110,5){\circle*{1.2}}
\put(115,5){\circle*{1.2}}
\put(107.5,10){\circle{1.2}}
\put(105,15){\circle*{1.2}}
\put(112.5,15){\circle*{1.2}}

\put(105,5){\line(1,2){2.3}}
\put(110,5){\line(-1,2){2.3}}
\put(107.25,10.5){\line(-1,2){2}}
\put(107.75,10.5){\line(1,1){4.5}}
\put(115,5){\line(-1,4){2.5}}

\put(105,2){\makebox(0,0)[cc]{$a$}}
\put(110,2){\makebox(0,0)[cc]{$a_1$}}
\put(115,2){\makebox(0,0)[cc]{$c$}}
\put(104,10){\makebox(0,0)[cc]{$x$}}
\put(105,18){\makebox(0,0)[cc]{$b$}}
\put(112.5,18){\makebox(0,0)[cc]{$b_1$}}

\put(125,5){\circle*{1.2}}
\put(130,5){\circle*{1.2}}
\put(135,5){\circle*{1.2}}
\put(130,10){\circle{1.2}}
\put(125,15){\circle*{1.2}}
\put(135,15){\circle*{1.2}}

\put(125,5){\line(1,1){4.5}}
\put(130,5){\line(0,1){4.5}}
\put(135,5){\line(-1,1){4.5}}
\put(129.5,10.5){\line(-1,1){4.5}}
\put(130.5,10.5){\line(1,1){4.5}}

\put(125,2){\makebox(0,0)[cc]{$a$}}
\put(130,2){\makebox(0,0)[cc]{$a_1$}}
\put(135,2){\makebox(0,0)[cc]{$c$}}
\put(126,10){\makebox(0,0)[cc]{$x$}}
\put(125,18){\makebox(0,0)[cc]{$b$}}
\put(135,18){\makebox(0,0)[cc]{$b_1$}}

\end{picture}

\caption{Proposition \ref{3.6}, (a) the case $|[a, b]| \ge 4$ (b) the case $|[a, b]| = 3$}
\end{figure}
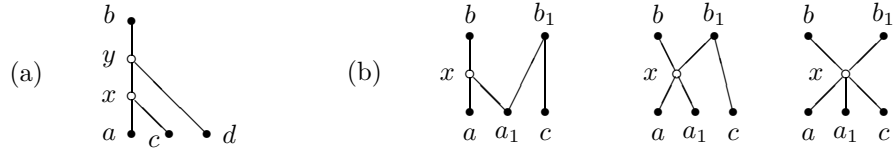

\vspace{.1in}

Consequently, the problem of classifying the locally 2-arc-transitive bipartite graphs with finite chain intervals is equivalent to that of 
classifying semilinear spaces satisfying the point-line-point, and line-point-line transitivity requirement given in the above proposition. As far 
as the authors are aware no such classification exists, although some subclasses of this family of semilinear spaces have been completely 
described. For example in \cite{Devillers1} the finite connected 4-homogeneous semilinear spaces are classified. Here a semilinear space is 
called \emph{$d$-homogeneous} if whenever the semilinear spaces induced on two subsets $S_1$ and  $S_2$ of $S$ of cardinality $\leq d$ are 
isomorphic there is an automorphism of $S$ mapping $S_1$ to $S_2$. (This is Devillers' terminology, which is a little different from our use of 
related terms.) If every such isomorphism extends we say $S$ is $d$-ultrahomogeneous. Clearly any 4-ultrahomogeneous semilinear space satisfies 
our point-line-point, and line-point-line transitivity condition, and so Devillers' classification gives rise to some examples. Another source of 
examples includes \cite{Devillers2}. See also \cite{Giudici3}. Some countably infinite examples of semilinear spaces satisfying the conditions of 
Proposition~\ref{3.6} arise from the work of K. Tent on `generalized $n$-gons'. A {\em generalized $n$-gon} is a bipartite graph such that the 
diameter of the graph is $n$ and there are no simple cycles (that is, without repetitions) of length less than $2n$. Also the graph is required to be 
thick, meaning that any element is incident with at least 3 other elements. An ordinary $n$-gon is a simple cycle of length $2n$ (in our 
poset-theoretic language we call these \emph{$2n$-crowns}, which are not however thick). In \cite{Tent}, using Hrushovski type constructions, Tent 
constructs for each $n \geq 3$ infinitely many non-isomorphic countable generalized $n$-gons for which the automorphism group acts transitively on the 
set of ordered ordinary $(n+1)$-gons contained in it. To obtain $2^{\aleph_0}$ examples one can consider projective planes over countable fields. 
There are clearly $2^{\aleph_0}$ pairwise non-isomorphic countable fields (for instance, of the form ${\mathbb Q}[\sqrt{p}: p \in P]$ for arbitrary 
sets $P$ of primes), and by the methods of von Staudt (see \cite{Veblen} Chapter 6) non-isomorphic fields give rise to non-isomorphic projective planes.

\subsection{Infinite chain diamond-free partial orders}
\label{subsec_chain}

We begin by identifying which chains may arise as intervals. In fact, here we find it convenient to work with $M^+$, and the corresponding interval 
$I^+(M)$, rather than $M^D$ and the interval $I(M)$. Once the possible $I^+(M)$ intervals have been determined, then by taking their completions we 
obtain the possible $I(M)$ intervals (note that this depends on the fact that we are in the chain interval case, see Lemma \ref{3.1}). 

We begin with a lemma about ramification points. The following result generalizes \cite[Lemma~2.3.11]{Warren} and is proved in a similar style to 
Lemma \ref{3.1}. 

\begin{lemma} \label{3.7}
Let $M$ be a partial order such that the intervals of $M^D$ are chains. Then 
\[
\uparrow \hspace{-.05in} \mathrm{Ram}(M) = \hspace{.05in}\uparrow \hspace{-.05in} \mathrm{Ram}(M^D) \quad \mbox{and} \quad 
\downarrow \hspace{-.05in} \mathrm{Ram}(M) = \hspace{.05in} \downarrow \hspace{-.05in} \mathrm{Ram}(M^D).
\]
\end{lemma}
\noindent{\bf Proof}: 
Since $M \subseteq M^D = (M^D)^D$ we have $\uparrow \hspace{-.05in} \mathrm{Ram}(M) \subseteq \hspace{.05in} \uparrow \hspace{-.05in} 
\mathrm{Ram}(M^D)$. 

For the converse let $a \in \hspace{.05in} \uparrow \hspace{-.05in} \mathrm{Ram}(M^D)$ be arbitrary, so $a = b \wedge c$ for $b, c \in M^D$ with $b 
\parallel c$. Let $x, y \in M$ with $x \geq b$ and $y \geq c$. It suffices to show that $x \wedge y = a$. 

First note that $b \not \le y$, since if $b \le y$ then $b, c \in [a, y]$, which is assumed to be a chain, and hence $b$ and $c$ are comparable, which 
is a contradiction. Similarly, $c \not \le x$. We deduce that $x \parallel y$. For suppose for instance that $x \le y$ (and $y \le x$ is similar). 
Then $x, c \in [a, y]$ which implies that $x$ and $c$ are comparable. Since $c \not \le x$, we can only have $x < c$, but then also $b < c$, a 
contradiction. It follows that $x \wedge y$ is an upward ramification point. Now $a \le b, x \wedge y \le x$, and as $[a, x]$ is linear, $b$ and 
$x \wedge y$ are comparable. Since $b \not \le y$, also $b \not \le x \wedge y$, and so $x \wedge y \le b$. Similarly, $x \wedge y \le c$. Hence $x 
\wedge y \le b \wedge c = a$, and so $x \wedge y = a$.    $\Box$

\vspace{.1in}

The above lemma relies heavily on the fact that the completion does not embed any diamonds, and it does not hold for general partial orders; see 
\cite[Section~4]{Gray1} for a discussion of this.  

We make use of Morel's classification of countable transitive linear orders. This involves the set $\mathbb{Z}^{\alpha}$, which is defined to be the 
family of all functions from $\alpha$ into $\mathbb Z$ with finite support, ordered lexicographically.

\begin{theorem}(Morel \cite{Morel})  \label{3.8}
A countable linear ordering is 1-transitive if and only if it is isomorphic to $\mathbb{Z}^{\alpha}$ or $\mathbb{Q}.\mathbb{Z}^{\alpha}$, for some 
countable ordinal $\alpha$, and where $\mathbb{Q}.\mathbb{Z}^{\alpha}$ is the lexicographic product of $\mathbb{Q}$ and $\mathbb{Z}^{\alpha}$. 
\end{theorem}

Following the notation used in \cite{Droste3}, we let $\mathbb{Q}_m$ be an $m$-coloured version of the rationals, as described in \cite{Warren}. This 
may be characterized as a structure of the form $(\mathbb{Q}, \leq, P_0, P_1, \ldots, P_{m-1})$ where $P_i$ are subsets which form a partition of 
$\mathbb{Q}$ into dense subsets, meaning that each non-empty open interval contains points of each $P_i$. We think of the sets $P_i$ as colours. Let 
$\mathbb{Q}_m(Z)$ be obtained by replacing each point of one fixed colour by the linear order $Z$. If $X$ and $Y$ are linear orders we 
write $X.Y$ to mean the `lexicographic product' $X$ copies of $Y$. We write $1 + X + 1$ for the linear order obtained from $X$ by adding two extra points, one at the top and one at the bottom.

We can derive the following result exactly as in \cite{Warren}. There the assumption was cycle-freeness, but actually only `diamond-freeness' is 
needed. Here we just give the result for the 3-$CS$-homogeneous case, which restricts us to the first three cases of the skeletal classification
\cite{Warren}, $\cal A$, $\cal B$, and $\cal C$. Similar results would hold for 3-$CS$-transitivity. 

\begin{lemma}\label{3.9}
Let $M$ be a countable $3$-$CS$-homogeneous $2$-level poset such that $I(M)$ is an infinite chain. Then $I^+(M) \setminus M$ is an infinite 
$2$-coloured chain (where the two `colours' correspond to $\uparrow \hspace{-.05in} \mathrm{Ram}(M)$ and $\downarrow \hspace{-.05in} \mathrm{Ram}(M)$) 
and exactly one of the the two following occurs:

(i) $\uparrow \hspace{-.05in} \mathrm{Ram}(M) = \hspace{.05in} \downarrow \hspace{-.05in} \mathrm{Ram}(M) = M^+ \setminus M$, and $I^+(M)$ is isomorphic to one of 
$1 + \mathbb{Z}^{\alpha} + 1$, or $1 + \mathbb{Q}.\mathbb{Z}^{\beta} + 1$, where $\alpha$, $\beta$ are countable ordinals such that $\alpha \geq 1$ 
(and all points which are not endpoints are coloured by both colours). 

(ii) $\uparrow \hspace{-.05in} \mathrm{Ram}(M) \hspace{.05in} \cap \downarrow \hspace{-.05in} \mathrm{Ram}(M) = \varnothing$, and $I^+(M)$ is 
isomorphic to one of the following:

(a) $1 + \mathbb{Q}.2 + 1$ (where $\mathbb{Q}.2$ is $\mathbb{Q}$ copies of the $2$-element chain $2$ with endpoints), where the lower point of each pair lies in $\uparrow 
\hspace{-.05in} \mathrm{Ram}(M)$ and the upper point of each pair lies in $\downarrow \hspace{-.05in} \mathrm{Ram}(M)$; 

(b) $1 + \mathbb{Q}_2 + 1$ (the $2$-coloured version of the rationals with endpoints).
\end{lemma}
\noindent{\bf Proof}: Since every point of $I(M)$ is a the least upper bound of a subset of $I^+(M)$, it follows that $I^+(M)$ is also infinite. The 
main point for the rest is that ${\rm Aut}(M)$ acts transitively on each of $\uparrow \hspace{-.05in} \mathrm{Ram}(M)$ and 
$\downarrow \hspace{-.05in} \mathrm{Ram}(M)$, as follows from 3-$CS$-transitivity, and also, that these two sets are countable. In a little more 
detail, the setwise stabilizer of $[a, b]$ in ${\rm Aut}(M)$ acts transitively on the family of its upward ramification points (similarly, its 
downward ramification points). For if $x, y \in \hspace{.05in} \uparrow \hspace{-.05in} \mathrm{Ram}(M)$, we may let $x = b \wedge c$ and 
$y = b \wedge d$, and by $3$-$CS$-homogeneity take $\{a, b, c\}$ to $\{a, b, d\}$, in that order, and then $[a, b]$ is fixed setwise, and $x$ must be 
mapped to $y$. We deduce that if $\uparrow \hspace{-.05in} \mathrm{Ram}(M)$ and $\downarrow \hspace{-.05in} \mathrm{Ram}(M)$ overlap at all, they are 
equal, which gives rise to the two cases described.

(i) In this case by the remark just made, the whole of $(a, b)$ forms a countable 1-transitive linear order and so it appears in the classification of Morel given above.

(ii) Now $\uparrow \hspace{-.05in} \mathrm{Ram}(M)$ and $\downarrow \hspace{-.05in} \mathrm{Ram}(M)$ are disjoint. Suppose first that for some $x \in 
\hspace{.05in} \uparrow \hspace{-.05in} \mathrm{Ram}(M)$ and $y \in \hspace{.05in} \downarrow \hspace{-.05in} \mathrm{Ram}(M)$, $x < y$ with no point 
in between. We call such a configuration a `pair'. Then by 3-$CS$-homogeneity, all ramification points lie in pairs. By Lemma \ref{3.2}, the family 
of pairs is densely linearly ordered without endpoints, and so the order-type of $I^+(M)$ is $\mathbb Q$ pairs with endpoints. Otherwise, between any 
two ramification points (of either kind, using Lemma \ref{3.2} again) there is another, of each type, and so the ordering is a copy of ${\mathbb Q}_2$ 
(with endpoints).    $\Box$

\vspace{.1in}

We shall now go on to consider each of the possible intervals $I$ identified in Lemma~\ref{3.9}, in each case investigating the class of countable locally $2$-arc-transitive bipartite graphs $M$ satisfying $I = I^+(M)$. This splits into three parts: (a) the case $I=1 + \mathbb{Z} + 1$ (b) more generally what we shall call the discrete interval types, which is where $I$ is isomorphic to $1 + \mathbb{Z}^{\alpha} + 1$, $1 + \mathbb{Q}.\mathbb{Z}^{\alpha} + 1$ (for some countable ordinal $\alpha \ge 1$), or to $1 + \mathbb{Q}.2 + 1$, and (c) the remaining dense intervals listed in Lemma~\ref{3.9}, namely the cases that $I$ is isomorphic to $1 + {\mathbb Q} + 1$ or to $1 + {\mathbb Q}_2 + 1$. 

Case (a) is considered first, and will lead to Theorem~\ref{3.15} where we show how the classification problem 
for that family relates to the problem of determining a certain family of highly arc-transitive digraphs. We then consider Case (b) in Theorem~\ref{3.16} where the direct connection with highly-arc-transitive digraphs no longer exists, but where the intuition developed in Case (a) allows us to construct many examples for each of these interval types; see Corollary~\ref{3.17}. Currently our methods do not extend to deal with the dense intervals of Case (c); see Remark~\ref{3.19} below for further discussion of this.

\subsection{\boldmath The case where the intervals are isomorphic to $1 + \mathbb{Z} + 1$} 
In this case the intervals in $M^+$ and in $M^D$ are the same since $1+ \mathbb{Z} + 1$ is already $DM$-complete. As we shall see below, in this case 
there is a close connection with highly arc-transitive digraphs, in the sense of \cite{Cameron2}. We now give the necessary background about these.

The definition of `digraph' was given at the beginning of Section~2, and we recall that the relation is required to be irreflexive and asymmetric, so 
we disallow loops, or arcs in both directions between any pair of vertices. The definitions of `$s$-arc', and `$s$-arc-transitive' for graphs carry 
over to digraphs and a digraph is said to be \emph{highly arc-transitive} if it is $s$-arc-transitive for all $s \ge 1$. For any vertex $v$ of a 
digraph $D$ we write $D^+(v) = \{ w : v \rightarrow w \}$ and $D^-(v) = \{ w : v \leftarrow w  \}$. 

\subsection{Reachability relations and descendants}

Given a 1-arc-transitive digraph $D$ one natural substructure that may be considered is that obtained by following an alternating walk. An 
\emph{alternating walk} in $D$ is a sequence $(x_0, \ldots, x_n)$ of vertices of $D$ such that either all $(x_{2i-1},x_{2i})$ and $(x_{2i+1},x_{2i})$ 
are arcs, or all $(x_{2i},x_{2i-1})$ and $(x_{2i},x_{2i+1})$ are arcs. If $a$ and $a'$ are arcs in $D$ and there is an alternating walk
$(x_0,\ldots,x_n)$ such that $(x_0,x_1)$ is $a$ and either $(x_{n-1},x_n)$ or $(x_n,x_{n-1})$ is $a'$, then $a'$ is said to be \emph{reachable} from 
$a$ by an alternating walk. This is denoted by $a \mathcal{A} a'$. Clearly the relation $\mathcal{A}$ is an equivalence relation on  $ED$, and the 
equivalence class containing the arc $a$ is denoted by $\mathcal{A}(a)$. If $D$ is $1$-arc-transitive then all digraphs of the form $\mathcal{A}(a)$ 
for $a \in ED$ are isomorphic, and we write this as $\Delta(D)$. The following result was given in \cite{Cameron2}.

\begin{prop}\cite[Proposition~1.1]{Cameron2} \label{3.10}
Let $D$ be a connected $1$-arc-transitive digraph. Then $\Delta(D)$ is $1$-arc-transitive and connected. Further, either

(i) $\mathcal{A}$ is the universal relation on $ED$ and $\Delta(D) = D$, or

(ii) $\Delta(D)$ is bipartite. 
\end{prop}

Another natural substructure of $D$ is given by the notion of `descendants'. For a vertex $u$ in $D$ a \emph{descendant} of $u$ is a vertex $v$ such 
that $D$ contains a directed path from $u$ to $v$. The set of all descendants of $u$ is denoted by $\mathrm{desc}(u)$. There is also the obvious dual 
notion of the \emph{ancestors} of a vertex. For $A \subseteq VD$ we define $\mathrm{desc}(A) = \bigcup_{v \in A}{\mathrm{desc}(v)}$. The set of 
\emph{ancestors} $\mathrm{anc}(v)$ of a vertex $v$ is the set of those vertices of $D$ for which $v$ is a descendant.

The $(m,n)$-directed trees discussed in Subsection~\ref{subsec_CFPO} are particularly simple examples of highly arc-transitive digraphs. We saw in that subsection how given such a directed tree $T$ countably many minimal points $X$ may be adjoined below certain maximal chains of $T$, and countably many maximal points $Y$ adjoined above, in such a way that $T$ may be recovered by taking the Dedekind--MacNeille completion of the two-level partial order induced by $M = X \cup Y$. Our first main result, Theorem~\ref{3.15}, will show that this is true in far greater generality, by describing a large class of countable highly arc-transitive digraphs for which it is possible to adjoin maximal and minimal points in this way.

\vspace{.1in}

Before proving that result, we shall need the following lemma concerning Dedekind--MacNeille completions in the case where the interval is of type $1 + \mathbb{Z} + 1$.

\begin{lemma}\label{3.12} Let $M$ be a countable locally $2$-arc-transitive bipartite graph such that the intervals of $M^+$ are isomorphic to 
$1 + \mathbb{Z} + 1$. Then 
\[
\uparrow \hspace{-.05in} \mathrm{Ram}(M) = \hspace{.05in} \downarrow \hspace{-.05in} \mathrm{Ram}(M) = M^+ \setminus M = M^D \setminus M.  
\]  
\end{lemma}
\noindent{\bf Proof}: The fact that $\uparrow \hspace{-.05in} \mathrm{Ram}(M), \downarrow \hspace{-.05in} \mathrm{Ram}(M) \subseteq M^+ \subseteq M^D$ 
is immediate. To see that $M \cap \downarrow \hspace{-.05in} \mathrm{Ram}(M)$ \newline $= \emptyset$, (and similarly, $M \cap 
\hspace{-.05in} \uparrow \hspace{-.05in} \mathrm{Ram}(M) = \emptyset$), let $a \in M \cap \hspace{-.05in} \uparrow \hspace{-.05in} \mathrm{Ram}(M)$. 
By applying Lemma \ref{3.2} to consecutive points of a copy of $\mathbb Z$ in $M^+$, we see that there is some upward ramification point $b$ which is 
not minimal. By 3-$CS$-transitivity of $M$ (which is the partial order version of local 2-arc-transitivity of the bipartite graph) there is an 
automorphism taking $a$ to $b$, which is clearly impossible.

Conversely, first suppose for a contradiction that $a \in M^D \setminus M^+$. Then $a$ is an ideal of $M$, so if we write $X$ and $Y$ for the sets of 
minimal and maximal elements of $M$, respectively, $a \subseteq X \cup Y$. If $a \cap Y \neq \emptyset$, then as $a^+ \neq \emptyset$, $|a \cap Y| = 
1$, and $a$ is principal, so lies in $M$, which is a contradiction. Hence $\emptyset \neq a \subseteq X$. If $|a| = 1$, then again $a \in M$, which is 
not the case. Hence $|a| \ge 2$, and applying a dual argument, $a^+$ is a subset of $Y$ of size at least 2. Pick distinct $x_0, x_1$ in $X$ and 
distinct $y_0, y_1$ in $Y$ such that $x_0, x_1 \le a \le y_0, y_1$. Thus $x_0 \vee x_1 \le a \le y_0 \wedge y_1$, and since $a \not \in M^+$, these 
inequalities are strict. Now $x_0 \vee x_1$ and $y_0 \wedge y_1$ lie in a copy of $\mathbb Z$ in $M^+$, and so there is a finite distance between them 
in $M^+$. Let $z \in M^+$ be the least member of this copy of $\mathbb Z$ above $a$. By Lemma \ref{3.2}, $[a, z) \cap (M \cup \uparrow \hspace{-.05in} 
{\rm Ram}(M))$ is non-empty, which contradicts the minimality of $z$. 

It remains to show that $\uparrow \hspace{-.05in} {\rm Ram}(M) = \hspace{.05in} \downarrow \hspace{-.05in} {\rm Ram}(M)$. Let 
$a \in \hspace{.05in} \downarrow \hspace{-.05in} {\rm Ram}(M)$. Then $a = x_0 \vee x_1$ for distinct $x_0, x_1 \in X$, and $a$ lies in a copy of 
$\mathbb Z$ in $M^+$. If we let $z$ be the successor of $a$ in this copy, then by applying Lemma \ref{3.2} to $[a, z)$ again, we find that 
$a \in \uparrow \hspace{-.05in} {\rm Ram}(M)$. Hence $\downarrow \hspace{-.05in} {\rm Ram}(M) \subseteq \uparrow \hspace{-.05in} {\rm Ram}(M)$, and 
the proof that $\uparrow \hspace{-.05in} {\rm Ram}(M) \subseteq \downarrow \hspace{-.05in} {\rm Ram}(M)$ is similar.    $\Box$

\vspace{.1in}

Note that this does not generalize to other intervals, like $\mathbb{Q}$ for example, since in those cases $M^+$ is countable, while $M^D \setminus M$ 
is not. Furthermore, we cannot drop the hypothesis of local 2-arc-transitivity, as we can `join together' two copies of a bipartite graph at just one 
maximal point, and this would violate the condition $M \cap \downarrow \hspace{-.05in} \mathrm{Ram}(M) = \emptyset$.

\begin{definition}[Intersection property] \label{3.13} 
We say that a digraph $D$ has the \emph{intersection property} if the intersection of any two principal ideals of $D$ is principal. More precisely, 
for all $x,y \in D$ if $\mathrm{desc}(x) \cap \mathrm{desc}(y) \neq \varnothing$ then there exists $z \in \mathrm{desc}(x) \cap \mathrm{desc}(y)$ such 
that $\mathrm{desc}(x) \cap \mathrm{desc}(y) = \mathrm{desc}(z)$.  
\end{definition}

\begin{definition}[Strongly transitive] \label{3.14}
A \emph{directed line} is a set $(x_i)_{i \in \mathbb{Z}}$ indexed by $\mathbb Z$ such that each $(x_i, x_{i+1})$ is a edge. 
By a \emph{$Y$-configuration} in a digraph $D$ we mean the digraph which is the amalgam of two infinite directed lines $(x_i)_{i \in \mathbb{Z}}$ and 
$(y_i)_{i \in \mathbb{Z}}$ amalgamated via the rule $x_i = y_i$ for all $i \geq 0$. Dually we have a $\overline{Y}$-configuration where we amalgamate 
via the rule $x_i = y_i$ for all $i \leq 0$. Then $D$ is said to be {\em strongly transitive} if its automorphism group acts transitively on the 
family of $Y$-configurations, and also on the family of $\overline{Y}$-configurations. 
\end{definition}

\begin{theorem}\label{3.15}
Let $D$ be a connected highly arc-transitive digraph such that

(i) the subdigraph induced on $\mathrm{desc}(u)$ is a tree and $2 \le {\rm outdegree}(u) \le \aleph_0$, for all $u \in VD$ (and the dual statement for 
$\mathrm{anc}(u)$);

(ii) $D$ has the intersection property;

(iii) $D$ is strongly transitive.

Then there is a countable connected $3$-$CS$-homogeneous $2$-level partial order $M$ such that $M^+ \setminus M \cong D$. 

Conversely if $\Gamma$ is a countable connected locally $2$-arc-transitive bipartite graph such that the interval $I(M)$ of $M = P(\Gamma)$ is 
isomorphic to $1 + {\mathbb Z} + 1$, then the digraph $D(M^+ \setminus M)$ naturally defined from the partial order $M^+ \setminus M$ is a connected 
highly arc-transitive digraph satisfying properties (i) and (ii).   
Furthermore, if $D(M^+ \setminus M)$ is locally finite then $D(M^+ \setminus M)$ satisfies (iii). 
\end{theorem}

\noindent{\bf Proof}: Let $D = (D, \to)$ be a digraph satisfying the stated conditions, and let $P = (P(D), \leq)$ be the corresponding partial order 
($P$ is a poset since anti-symmetry is automatically satisfied by (i)). It follows from the assumptions that $P$ does not have any maximal or minimal 
elements. Our aim is to extend $P$ to a countable poset $\widehat{P} = X \cup P \cup Y$ with the following properties:
\begin{itemize}
\item $X = \mathrm{Min}(\widehat{P})$, $Y = \mathrm{Max}(\widehat{P})$;
\item $(\forall p \in P)(\exists x \in X)(\exists y \in Y)(x \leq p \leq y)$; 
\item $\widehat{P} \cong M^D \cong M^+$ where $M$ is the substructure of $\widehat{P}$ induced on $X \cup Y$;  
\item $\Gamma(M)$ is a connected countable locally 2-arc transitive bipartite graph.
\end{itemize}
We achieve this in two stages. First we build an extension $E = A \cup P \cup B$ of $P$ satisfying all the desired properties, except that $E$ will be 
uncountable. Then we show how to cut down to a countable substructure $\widehat{P}$ with $P \leq \widehat{P} \leq E$ and such that $\widehat{P}$ still 
satisfies all the desired properties. This is very much in the spirit of the downward L\"{o}wenheim--Skolem Theorem; see \cite{Marker}.

\vspace{.1in}

\noindent{\bf Stage 1}: Building an uncountable extension $E$:

The idea is to add points above and below $D$ corresponding to all possible sequences. To do this formally, define a \emph{ray} in $D$ to be an 
infinite sequence $(v_i: i \in \mathbb{N})$ of vertices such that $v_i \rightarrow v_{i+1}$ for all $i$. Dually we define an \emph{antiray} as a 
sequence $(v_i : i \in \mathbb{Z}^{ < 0})$ such that $v_i \rightarrow v_{i+1}$ for all $i \leq -1$. Now we want to add new minimal points below each 
ray of $P = P(D)$, but distinct rays may lie above the same minimal point, depending on whether or not they `eventually agree', so we need to define a 
relation $\sim$ on the set of all rays by letting $R_1 \sim R_2$ where $R_1 = (r_i : i \in \mathbb{N})$ and $R_2 = (s_i : i \in \mathbb{N})$ if for 
some $n, m \in \mathbb{N}$ and all $j \ge 0$, $r_{n+j} = s_{m+j}$ (i.e. $R_1$ and $R_2$ `eventually agree'). Clearly $\sim$ is an equivalence relation 
on the set of all rays, and we let $A$ denote the set of all $\sim$-classes of rays of $P$, and extend $\leq$ from $P$ to $P \cup A$ by defining $a < 
p$ for $a \in A$ and $p \in P$ if $p \in R$ for some ray $R \in a$. There is an obvious dual equivalence relation $\approx$ on the collection of all 
antirays; we let $B$ be the set of $\approx$-classes of antirays, and let $b > p$ for $b \in B$ and $p \in P$ if $p \in R$ for some antiray $R \in b$.
Finally for $a \in A$ and $b \in B$ we let $a \leq b$ if $a \leq p \leq b$ for some $p \in P$. 

It is easy to see that, defined in this way, $E = (A \cup P \cup B, \leq)$ is a poset with $\mathrm{Max}(E) = B$, $\mathrm{Min}(E) = A$, and that the 
original poset $(P,\leq)$ naturally embeds in $E$ as the substructure induced by $P$. 

Let $N$ be the substructure of $E$ induced on $A \cup B$. Then $N$ is a poset with maximal chains of height 2, and as such, may naturally be viewed as 
a bipartite graph $\Gamma(N)$. In fact we are now able to show that $\Gamma(N)$ is an uncountable connected locally 2-arc transitive bipartite graph.
Let $v$ be a vertex of $D$. By the assumption that $\mathrm{desc}(v)$ is a tree, whose vertices have outdegree at least 2 but no greater than 
$\aleph_0$, it follows that there are $2^{\aleph_0}$ distinct rays emanating from $v$, belonging to pairwise distinct $\sim$-classes, which shows that 
$|A| = 2^{\aleph_0}$, and similarly, $|B| = 2^{\aleph_0}$. Connectedness of $\Gamma(N)$ follows from that of $D$. Finally, we must show that 
$\Gamma(N)$ is locally 2-arc transitive. Clearly every automorphism $\alpha \in \mathrm{Aut}(D) = \mathrm{Aut}(P(D))$ extends to an automorphism 
$\bar{\alpha} \in \mathrm{Aut}(E) = \mathrm{Aut}(A \cup P \cup B)$, which in turn induces an automorphism of $N = A \cup B$ and hence of $\Gamma(N)$. 
Let $\gamma: (u,v,w) \mapsto (u',v',w')$ be an isomorphism between substructures of $\Gamma(N)$, where $u,w,u',w' \in B$ and $v,v' \in A$, $v$ 
adjacent to $u$ and $w$, and $v'$ adjacent to $u'$ and $w'$. Choose antiray and ray representatives $R_u, \ldots, R_{w'}$ of the respective 
$\approx$- and $\sim$-classes. It is an easy consequence of the assumptions that the antirays $R_u$, $R_w$ and the ray $R_v$ may be chosen to meet at 
just a single point, and in such a way that $Y = R_u \cup R_v \cup R_w$ is a $Y$-configuration. Similarly we may suppose that 
$Y' = R_u' \cup R_v' \cup R_w'$ is a $Y$-configuration. By assumption, $D$ is transitive on $Y$-configurations and thus there is an automorphism 
$\alpha$ of $D$ such that $\alpha: Y \rightarrow Y'$, and then the automorphism of $\Gamma(N)$ that $\alpha$ induces extends the isomorphism $\gamma$. 
Along with a dual argument, applying $\overline{Y}$-transitivity, this completes the proof that $\Gamma(N)$ is locally 2-arc transitive.    

\vspace{.1in}

\noindent{\bf Stage 2}: Cutting down to a countable structure:

Our aim is to construct a poset $M = X \cup Y$, where $X$ and $Y$ are countable subsets of $A$ and $B$ respectively, and $M$ satisfies all the 
desired properties given in the statement of the theorem. The sets $M, X$ and $Y$ will be defined in countably many steps. At the same time we shall 
define a subgroup $G$ of $\mathrm{Aut}(N)$ that will act locally 2-arc transitively on $M$. 

For each of the countably many vertices $p \in P = P(D)$, fix $p_{\mathrm{max}} \in B$ and $p_{\mathrm{min}} \in A$ such that $p_{\mathrm{min}} < p < 
p_{\mathrm{max}}$. Then let \[ X_0 = \{ p_{\mathrm{min}} : p \in P \}, \quad Y_0 = \{ p_{\mathrm{max}} : p \in P \}, \quad M_0 = X_0 \cup Y_0.
\]
Note that $M_0$ is countable since $P$ is countable. Enumerate all isomorphisms between $\vee$-2-arcs of $M_0$. For each such choose and fix an 
automorphism of $N$ that extends the given isomorphism. This is possible since, as shown above, $N$ is locally 2-arc transitive. Do the same for 
the countably many isomorphisms between $\wedge$-2-arcs, and let $C_1 \subseteq \mathrm{Aut}(N)$ be the resulting countable set of automorphisms. Let 
$G_1$ be the subgroup of $\mathrm{Aut}(N)$ generated by $C_1$, and let 
\[
X_1 = X_0^{G_1} = \{ x^g : x \in X_0, g \in G_1 \}, \quad Y_1 = Y_0^{G_1},
\]      
noting that $X_0 \subseteq X_1 \subseteq A$, $Y_0 \subseteq Y_1 \subseteq B$ and $X_1$ and $Y_1$ are both countable. Let $M_1$ be the substructure of 
$N$ induced by $X_1 \cup Y_1$. 

At a typical stage of the construction we are given $M_i = X_i \cup Y_i$ and $G_i \leq \mathrm{Aut}(N)$. We enumerate all the countably many 
isomorphisms between $\vee$-2-arcs and $\wedge$-2-arcs, let $C_{i+1}$ be a countable set of automorphisms extending these isomorphisms, and let 
$G_{i+1} = \langle G_i \cup C_{i+1} \rangle \leq \mathrm{Aut}(N)$. Note that $G_{i+1}$ is countable and $G_i \leq G_{i+1}$. Then define $X_{i+1} = 
X_i^{G_{i+1}}$, $Y_{i+1} = Y_i^{G_{i+1}}$ and $M_{i+1}$ to be the substructure of $N$ induced by $X_{i+1} \cup Y_{i+1}$. Finally define:
\[
M = \bigcup_{i \geq 0} M_i, \quad G = \bigcup_{i \geq 0} G_i \leq \mathrm{Aut}(N). 
\]
Now let $g \in G$ and $m \in M$ be arbitrary, say with $m \in X$. Then there is some $i \in \mathbb{N}$ such that $g \in G_{i+1}$ and $m \in M_i$. Then from the definition:
\[
m^g \in X_i^g \subseteq X_i^{G_{i+1}} = X_{i+1} \subseteq X.
\]
It follows that $G \leq \mathrm{Aut}(N)$ fixes $M$ setwise, and hence $G$ acts on $M$ as a group of automorphisms. It now follows from the way $X$, 
$Y$ and $G$ were defined, that $G$ acts locally 2-arc transitively on $M$. 

Since $M_0 \subseteq M$ was chosen densely, and since every vertex in $D$ has in- and out-degree strictly greater than one, it follows that $P = P(D) 
\subseteq M^D$. Thus we are left with the task of showing that $M \cup P$ is $DM$-complete. Observe that the maximal chains have order-type 
$\mathbb{Z}$, $1+ \mathbb{Z}$, $\mathbb{Z} + 1$ or $1+ \mathbb{Z} + 1$, and these are all $DM$-complete. Next consider $x, y \in M \cup P$ with $x 
\parallel y$ and such that they are bounded below in $M \cup P$. Suppose that $x,y \in Y = \mathrm{Max}(M)$ and that they are bounded below by $z \in 
X = \mathrm{Min}(M)$. Let $R_x$ and $R_y$ be antirays with $R_x \in x$ and $R_y \in y$, and let $R_z$ be a ray in $z$, chosen so that the sets $R_x$, 
$R_y$ and $R_z$ are disjoint and the union $R_x \cup R_y \cup R_z$ is a $Y$-configuration. It is then straightforward, using the intersection 
property, to verify that the unique maximal element of the ray $R_y$ is the unique greatest lower bound in $P \cup M$ of $x$ and $y$. There are 
several other cases, depending which of the sets $X, Y$ or $P$ the points $x$ and $y$ belong to. Each case may be dealt with using a similar argument 
to that used for the case considered above, and as a result we deduce that $M \cup P$ is indeed $DM$-complete. 

Conversely let $\Gamma$ be a countable connected locally 2-arc-transitive bipartite graph such that, with $M = P(\Gamma)$, the interval $I(M)$ of $M$ 
is isomorphic to $1 + \mathbb{Z} + 1$. Now consider the digraph $D = D(M^+ \setminus M)$, whose edges are pairs $(x, y)$ such that $x > y$ with no 
points in between. By Lemma \ref{3.12}, $M^+ \setminus M = M^D \setminus M$ and it follows from this and the assumption 
$I(M) \cong 1 + \mathbb{Z} + 1$ that in $D$ each of $\mathrm{desc}(u)$ and $\mathrm{anc}(u)$ are trees, and the cardinality restrictions are satisfied 
since $\Gamma$ is countable. Hence property (i) is satisfied. Property (ii), the intersection property, holds in $D$ since $M^D = M^+$ is 
$DM$-complete. It follows from local 2-arc-transitivity first of all that the automorphism group of $D$ acts transitively, since any point is a 
ramification point, and it acts highly arc-transitively for essentially the same reason, since any (directed) $s$-arc for finite $s$ is part of one of 
the intervals determined by a 2-arc of $\Gamma$. 

For the last clause, suppose that $D = D(M^+ \setminus M)$ is a locally finite digraph. The automorphism group of an infinite locally finite digraph is a topological group, with the topology of pointwise convergence. Since $D$ is locally finite, under this topology, for all $v \in VD$, the stabilizer $\mathrm{Aut}(D)_v$ is  compact and $\mathrm{Aut}(D)$ is locally compact; see for instance \cite[Section 3]{Evans}. Let $\phi: Y \rightarrow Y'$ be an isomorphism between two $Y$-configurations in the digraph $D$. We show that $\phi$ can be extended to an automorphism of $D$ (the argument for $\overline{Y}$-configurations is dual). Let $v \in Y$ be the unique vertex of the $Y$-configuration with out-degree 2. Write $Y = \bigcup_{i \geq 0} Y_i$ as an infinite union of finite $Y$-configurations (each containing $v$) with $Y_i \subset Y_{i+1}$ for all $i$. For each $i$ let $\phi_i :Y_i \rightarrow Y_i'$ be the restriction of the isomorphism $\phi$, where $Y_i' \subseteq Y'$ denotes the image of $Y_i$ under $\phi$. Let $v' = \phi(v)$. Since $\Gamma$ is locally $2$-arc-transitive by assumption, and since every automorphism of $\Gamma$ naturally induces an automorphism of the digraph $D$, we see that each isomorphism $\phi_i$ extends to an automorphism $\hat{\phi_i}$ of $D$ induced by an automorphism of $\Gamma$. Since each $\hat{\phi_i}$ satisfies $\hat{\phi_i}(v) = v'$ it follows that $\{ \hat{\phi_i} : i \geq 0 \}$ is a subset of $\mathrm{Aut}(D)_v \hat{\phi_0}$, where $\mathrm{Aut}(D)_v \hat{\phi_0}$ is compact since it is a coset of $\mathrm{Aut}(D)_v$. Since $\{ \hat{\phi_i} : i \geq 0 \}$ is infinite it has at least one accumulation point $\hat{\phi}$ in $\mathrm{Aut}(D)_v \hat{\phi_0}$. Then $\hat{\phi}$ lies in  $\mathrm{Aut}(D)_v \hat{\phi_0} \subseteq \mathrm{Aut}(D)$ and is an automorphism which extends the original isomorphism $\phi$. This proves that $\mathrm{Aut}(D)$ is transitive on $Y$-configurations. Taken together with the dual argument for $\overline{Y}$-configurations, this proves that $D$ satisfies (iii).    

Finally, we note that it is not clear whether we can prove strong transitivity of $D$ in all cases 
(which we would have liked to deduce from local 2-arc-transitivity of $\Gamma$) because not all Y-shapes have points of $M$ above and below their 
maximal chains (but fortunately (iii) is not needed in what follows).$\Box$

\vspace{.1in}

It is natural to ask whether there are actually any highly arc-transitive digraphs satisfying conditions (i), (ii) and (iii) of Theorem~\ref{3.15}. Any $(m,n)$-directed tree certainly satisfies these three conditions, and in these cases the result of the theorem amounts to the explanation of the $CFPO$ construction for these cases given in Subsection~\ref{subsec_CFPO}. In fact, many more examples may be constructed by utilizing the construction of universal highly arc-transitive digraphs from \cite{Cameron2}. We give a brief description of this construction here, referring the reader to \cite[Section~2]{Cameron2} for full details.

Let $\Delta$ be a 1-arc-transitive, connected bipartite graph with given bipartition $X \cup Y$. Let $u = |X|$ and $v = |Y|$, noting that $u$ and $v$ 
need not be finite. We shall construct a digraph $DL(\Delta)$ that has the property that its reachability graph is isomorphic to $\Delta$. Let $T$ be 
a directed tree with constant in-valency $v$ and constant out-valency $u$. For each vertex $t \in T$ let $\phi_t$ be a bijection from $T^{-}(t)$ to 
$Y$, and let $\psi_t$ be a bijection from $T^+(t)$ to $X$. Then $DL(\Delta)$ is defined to be the digraph with vertex set $ET$ such that for $(a,b), 
(c,d) \in ET$, $((a,b),(c,d))$ is a directed edge of $DL(\Delta)$ if and only if $b=c$ and $(\psi_b(a),\phi_b(d))$ is an edge of $\Delta$.  The graph 
$DL(\Delta)$ may be thought of as being constructed by taking $T$ and replacing each vertex of $T$ by a copy of $\Delta$. Then for copies of $\Delta$ 
that are indexed by adjacent vertices $a$ and $b$ of $T$ we identify a single vertex from one of the copies of $\Delta$ with a vertex from the other 
copy of $\Delta$ (with the bijections determining the identifications). 

Since $\Delta$ is 1-arc-transitive it follows that different choices of the bijections $\phi_y$ and $\psi_y$ for $y \in VT$ will lead to isomorphic 
digraphs, and so the definition of $DL(\Delta)$ is unambiguous.

\begin{figure}
\unitlength=1mm
\begin{picture}(120.00,70.00)

\put(50,71){\line(0,1){1}}
\put(50,73){\line(0,1){1}}
\put(50,75){\line(0,1){1}}

\put(90,71){\line(0,1){1}}
\put(90,73){\line(0,1){1}}
\put(90,75){\line(0,1){1}}

\put(70,61){\line(0,1){1}}
\put(70,63){\line(0,1){1}}
\put(70,65){\line(0,1){1}}

\put(90,47){\line(0,1){1}}
\put(90,49){\line(0,1){1}}
\put(90,51){\line(0,1){1}}

\put(110,41){\line(0,1){1}}
\put(110,43){\line(0,1){1}}
\put(110,45){\line(0,1){1}}

\put(110,17){\line(0,1){1}}
\put(110,19){\line(0,1){1}}
\put(110,21){\line(0,1){1}}

\put(50,25){\line(0,1){1}}
\put(50,27){\line(0,1){1}}
\put(50,29){\line(0,1){1}}

\put(70,7){\line(0,1){1}}
\put(70,9){\line(0,1){1}}
\put(70,11){\line(0,1){1}}

\put(40,65){\circle*{1.2}}
\put(50,65){\circle*{1.2}}
\put(60,65){\circle*{1.2}}
\put(80,65){\circle*{1.2}}
\put(90,65){\circle*{1.2}}
\put(100,65){\circle*{1.2}}
\put(40,55){\circle*{1.2}}
\put(50,55){\circle*{1.2}}
\put(60,55){\circle*{1.2}}
\put(70,55){\circle*{1.2}}
\put(80,55){\circle*{1.2}}
\put(90,55){\circle*{1.2}}
\put(100,55){\circle*{1.2}}
\put(40,45){\circle*{1.2}}
\put(50,45){\circle*{1.2}}
\put(60,45){\circle*{1.2}}
\put(70,45){\circle*{1.2}}
\put(80,45){\circle*{1.2}}
\put(90,45){\circle*{1.2}}
\put(100,45){\circle*{1.2}}
\put(40,35){\circle*{1.2}}
\put(50,35){\circle*{1.2}}
\put(60,35){\circle*{1.2}}
\put(80,35){\circle*{1.2}}
\put(90,35){\circle*{1.2}}
\put(100,35){\circle*{1.2}}
\put(110,35){\circle*{1.2}}
\put(120,35){\circle*{1.2}}
\put(100,25){\circle*{1.2}}
\put(110,25){\circle*{1.2}}
\put(120,25){\circle*{1.2}}

\put(65,15){\circle*{1.2}}
\put(70,15){\circle*{1.2}}
\put(75,15){\circle*{1.2}}
\put(65,20){\circle*{1.2}}
\put(75,20){\circle*{1.2}}

\put(40,55){\line(0,1){10}}
\put(40,55){\line(1,1){10}}
\put(50,55){\line(0,1){10}}
\put(50,55){\line(1,1){10}}
\put(60,55){\line(0,1){10}}
\put(60,55){\line(-2,1){20}}

\put(80,55){\line(0,1){10}}
\put(80,55){\line(1,1){10}}
\put(90,55){\line(0,1){10}}
\put(90,55){\line(1,1){10}}
\put(100,55){\line(0,1){10}}
\put(100,55){\line(-2,1){20}}

\put(60,45){\line(0,1){10}}
\put(60,45){\line(1,1){10}}
\put(70,45){\line(0,1){10}}
\put(70,45){\line(1,1){10}}
\put(80,45){\line(0,1){10}}
\put(80,45){\line(-2,1){20}}

\put(40,35){\line(0,1){10}}
\put(40,35){\line(1,1){10}}
\put(50,35){\line(0,1){10}}
\put(50,35){\line(1,1){10}}
\put(60,35){\line(0,1){10}}
\put(60,35){\line(-2,1){20}}

\put(80,35){\line(0,1){10}}
\put(80,35){\line(1,1){10}}
\put(90,35){\line(0,1){10}}
\put(90,35){\line(1,1){10}}
\put(100,35){\line(0,1){10}}
\put(100,35){\line(-2,1){20}}

\put(100,25){\line(0,1){10}}
\put(100,25){\line(1,1){10}}
\put(110,25){\line(0,1){10}}
\put(110,25){\line(1,1){10}}
\put(120,25){\line(0,1){10}}
\put(120,25){\line(-2,1){20}}

\put(65,15){\line(0,1){5}}
\put(65,15){\line(1,6){5}}
\put(70,15){\line(0,1){30}}
\put(70,15){\line(1,1){5}}
\put(75,15){\line(0,1){5}}
\put(75,15){\line(-2,1){10}}

\end{picture}

\caption{The construction using 6-cycles}\label{fig_6cycle}
\end{figure}
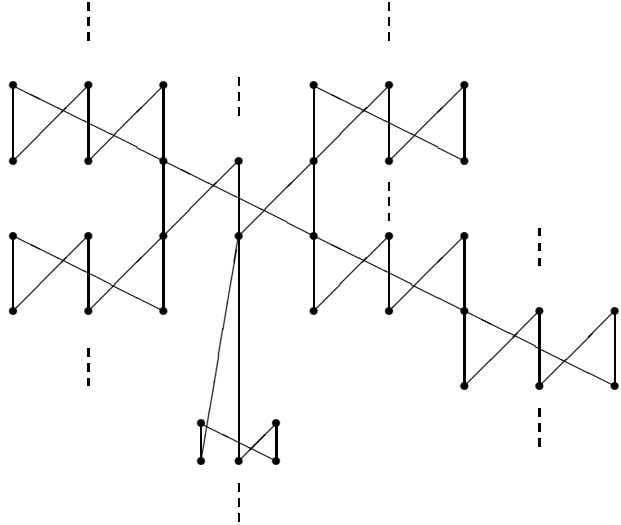

An illustration of $DL(\Delta)$ in the case that $\Delta$ is a bipartite $6$-cycle is given in Figure~\ref{fig_6cycle}. 
Note that $\Delta$ may be recovered from $DL(\Delta)$ via the reachability relation $\mathcal{A}$ discussed at the beginning of this subsection.
It is an easy consequence of the definition of $DL(\Delta)$ that for any $DM$-complete countable locally $2$-arc-transitive bipartite graph $\Delta$ (equivalently, the incidence graph of a semilinear space satisfying the transitivity property given in Proposition~\ref{3.6}), the digraph $DL(\Delta)$ is a highly arc-transitive digraph satisfying conditions (i), (ii) and (iii) of Theorem~\ref{3.15}. We saw in the discussion at the end of Subsection~\ref{subsec_finitechain} that there are continuum many different choices for such $\Delta$, each of which gives rise to a different digraph $DL(\Delta)$ which in each case, applying Theorem~\ref{3.15}, gives rise to a distinct $3$-$CS$-homogeneous $2$-level partial order $M$ such that $M^+ \setminus M \cong DL(\Delta)$. 

We shall now move on to discuss the problem for the other discrete intervals given in Lemma~\ref{3.9} and will generalize this statement in Corollary~\ref{3.17} showing that in each case continuum many pairwise non-isomorphic examples may be constructed.

\subsection{Axiomatic approach: the cases where the intervals are isomorphic to 
 $1 + \mathbb{Z}^{\alpha} + 1$, $1 + \mathbb{Q}.\mathbb{Z}^{\alpha} + 1$, or $1 + \mathbb{Q}.2 + 1$
} 

Now we describe an axiomatic approach to the construction described in Theorem \ref{3.15}, which will throw light on the general method and enable us 
to establish uniqueness of the constructed poset, given two `inputs', a $DM$-complete locally $2$-arc-transitive bipartite graph $\Delta$ and the linear order $\cal Z$. In addition, we 
generalize to all the cases listed in Lemma \ref{3.9} with the exception of the dense cases $\mathbb Q$ and 
${\mathbb Q}_2$, even though in many of these, the direct links with digraphs as explained in Theorem \ref{3.15} are more tenuous, and we have to move 
over to the partial order approach. Specifically, the digraphs will usually be disconnected, and the relationship between the connected components can 
only be recognized by using the partial order, and not just the digraph. In addition, we restrict to the cases in which the interval of $M$ has 
order-type $1 + {\cal Z} + 1$ where $\cal Z$ is ${\mathbb Z}^\alpha$, ${\mathbb Q}.{\mathbb Z}^\alpha$ for a countable ordinal $\alpha > 0$, or 
${\mathbb Q}.2$. The dense cases, $\mathbb Q$ and ${\mathbb Q}_2$, are more involved and we only make some remarks about these later describing some 
difficulties in extending our results to these cases. In the first list we have all those order-types which are not dense, and so we can envisage 
adjacent vertices forming part of a copy of a bipartite graph. Let us therefore fix a countable connected $DM$-complete locally 1- and 
2-arc-transitive bipartite graph $\Delta$ (for instance a 6-crown, as in Figure 4, though the method applies more generally than that, as we now see; 
by Proposition \ref{3.5}, $\Delta$ arises from a semilinear space). We remark that local 2-arc-transitivity of $\Delta$ implies that 
$\uparrow \hspace{-.05in} \mathrm{ro}(x)$ and $\downarrow \hspace{-.05in} \mathrm{ro}(x)$ are well-defined for minimal and maximal members of $\Delta$ 
respectively, and we write them as $\uparrow \hspace{-.05in} \mathrm{ro}(\Delta)$ and $\downarrow \hspace{-.05in} \mathrm{ro}(\Delta)$. In what 
follows when we talk of a `copy' of $\Delta$ in $Z = P \setminus (X \cup Y)$ we understand that it is embedded as a digraph (thus, strictly speaking, 
in $D(Z)$, in the notation of Theorem~\ref{3.15}), which is the same as saying in poset language, that consecutive points of the copy of $\Delta$ are also consecutive in $P$. 

Let $P$ be a connected countably infinite poset with $X = \mathrm{Min}(P)$, $Y = \mathrm{Max}(P)$ disjoint, and $Z = P \setminus (X \cup Y)$, 
let $\cal Z$ be a linear order equal to $\mathbb{Z}^{\alpha}$ 
or $\mathbb{Q}.\mathbb{Z}^{\alpha}$ for some countable ordinal $\alpha \ge 1$, or $\mathbb{Q}.2$,
and 
suppose that $P$ has the following properties:

(P1) $Z = {\rm Ram}(X \cup Y)$, and if ${\cal Z} = {\mathbb Q}.2$ then 
$\uparrow \hspace{-.05in} \mathrm{Ram}(Z) \hspace{.05in} \cap \downarrow \hspace{-.05in} \mathrm{Ram}(Z) = \emptyset$, 

(P2) $Z$ does not embed any non-alternating cycles,

(P3) any cycle that $Z$ embeds lies in a copy of $\Delta$,

(P4) for all $z \in Z$, if there is a point of $Z$ immediately above $z$, then there is a copy of $\Delta$ containing $z$ and all members of 
$z^\uparrow$ consecutive with $z$, and if there is a point of $Z$ immediately below $z$, then there is a copy of $\Delta$ containing $z$ and all 
members of $z^\downarrow$ consecutive with $z$ (in the case of finite ramification orders, this may be replaced by demanding that $\uparrow 
\hspace{-.05in} \mathrm{ro}(z) = \hspace{.05in} \uparrow \hspace{-.05in} \mathrm{ro}(\Delta)$ and $\downarrow \hspace{-.05in} \mathrm{ro}(z) = 
\hspace{.05in} \downarrow  \hspace{-.05in} \mathrm{ro}(\Delta)$),

(P5) every pair of consecutive vertices of $P$ belongs to a unique copy of $\Delta$,

(P6) for all $z \in Z$ there are $x \in X$ and $y \in Y$ such that $x < z < y$,

(P7) for any $x \in X$ and $y \in Y$ with $x < y$, $[x, y] \cong 1 + {\cal Z} + 1$, (so $P$ is diamond-free),

(P8) no maximal chain of $Z$ of the form $(x, y)$ for $x \in X$ and $y \in Y$ with $x < y$ has more than one point of $Y$ above it, or more than one 
point of $X$ below it.

\begin{theorem}\label{3.16}
For any countable connected $DM$-complete locally $2$-arc-transitive bipartite graph $\Delta$ and linear order $\cal Z$ equal to $\mathbb{Z}^{\alpha}$ 
or $\mathbb{Q}.\mathbb{Z}^{\alpha}$ for some countable ordinal $\alpha \ge 1$, or $\mathbb{Q}.2$, there is a partially ordered set $P$ which fulfils 
conditions (P1)-(P8), and $P$ is uniquely determined up to isomorphism. Furthermore, the $2$-level partial order $M = X \cup Y = M(\Delta, {\cal Z})$ 
is connected, countable, and $3$-$CS$-homogeneous, and if two partial orders arising from this construction are isomorphic, then the corresponding 
values of $\cal Z$ are isomorphic, and so are the corresponding values of $\Delta$. 
\end{theorem}

\noindent{\bf Proof}: To establish existence in the $1 + {\mathbb Z} + 1$ case, let $\Gamma = DL(\Delta)$, and let $P = 
P(\Gamma)$ be the partial order given by applying Theorem \ref{3.15}. In the general case, this has to be adapted, and the easiest way to build $P$ is 
via `approximations', which are needed in the uniqueness proof as well. By an {\em approximation} we understand an infinite connected diamond-free 
partially ordered set $(A, <)$, which is a finite union of maximal chains, all isomorphic to $1 + {\cal Z} + 1$, and copies of $\Delta$, such that if 
$x < y$ in such a copy, then $x$ and $y$ are consecutive in $A$, and such that if $z$ is an upward ramification point of $A$ then $z$ is a minimal 
point of a copy of $\Delta$ in $A$, and if it is a downward ramification point of $A$ then it is a maximal point of a copy of $\Delta$ in $A$. In 
addition, we require that properties (P2) and (P3) hold for $A$, and also that if $C$ is a maximal chain of $A$ in type $1 + {\cal Z} + 1$, then its 
greatest and least points $y$ and $x$ are the only points above and below $C \setminus \{x, y\}$ respectively (approximating property (P8)). To make 
things easier to handle, we also assume that an `approximation' comes with a homomorphism $l$ from $A$ onto $1 + {\cal Z} + 1$, which is meant to keep 
track of which levels the elements of $A$ lie on, and which is required to restrict to an isomorphism from any maximal chain of $A$ to 
$1 + {\cal Z} + 1$. We can start with an approximation just consisting of one copy of $1 + {\cal Z} + 1$, which clearly satisfies all these 
stipulations, where the value of $l$ is immediate, and the main point is to see how to extend to make all the properties true.

The main properties that we have to make true are (P4), (P5), and (P6). Each of these only involves countably many tasks, so provided we can do each 
on its own, then we may repeat (dovetailing in new tasks as they arise) and ensure that all possible tasks are fulfilled. We remark that existence in 
(P5) follows from (P4), and we just need to ensure uniqueness. For this, we just make sure that once a consecutive pair $z_1 < z_2$ in $A$ has been 
included in a copy of $\Delta$, we never add another copy of $\Delta$ containing these two points (in fact, we never add another copy of $\Delta$ 
having $z_1$ as a lower point or having $z_2$ as an upper point). For (P4) (and (P5)) the case we need to handle is therefore that in which 
$z_1 < z_2$ are consecutive members of $A$ such that $z_1$ is not the lower point of a copy of $\Delta$ and $z_2$ is not an upper point of a copy of 
$\Delta$, in $A$. We choose $u_1 < u_2$ in $\Delta$, and let $A'$ be the union of $A$ and $\Delta$ which is disjoint apart from identifying $z_1$ with 
$u_1$ and $z_2$ with $u_2$. This is partially ordered by the transitive closure of the union of the relations on $A$ and $\Delta$, and it is then 
clear that $A'$ is an approximation, when we assign the value $l(z_1)$ to all lower points of this copy of $\Delta$, and $l(z_2)$ to all its upper 
points. 

To extend so that (P6) holds, we suppose that $z \in A$ is given. If $x$ and $y$ as in (P6) already exist in $A$, then no extension is necessary. 
Otherwise if $x$ does not exist for instance, this implies that there is a minimal point $u$ of $A$ below $z$ (which must be a lower point of a copy 
of $\Delta$). In forming $A'$ this time, we adjoin points below $u$ in order-type $1 + {\cal Z}_1$ (where ${\cal Z}_1 = \{t \in {\cal Z}: t < 
l(u)\}$). If $y$ does not exist, then we also must consider a maximal point $v$ of $A$ above $z$ and add points above $v$ in order-type 
${\cal Z}_2 = \{t \in {\cal Z}: t > l(v)\}$. Here the choice of ${\cal Z}_1$ and ${\cal Z}_2$ guarantees that the resulting maximal chain has 
order-type $\cal Z$, and the choice of the extension of $l$ to the new points to keep it a homomorphism is immediate.

To establish uniqueness, let $P_1$ and $P_2$ be countable posets satisfying the properties (P1)-(P8). We want to prove that $P_1$ and $P_2$ are 
isomorphic. The isomorphism is built by back-and-forth in countably many steps. We say that a substructure $A$ of $P_i$ is {\em convex} if whenever $x 
< z < y$ and $x, y \in A$, then also $z \in A$. The notion of `approximation' is nearly the same as before, except that this time it will be a 
substructure of $P_1$ or $P_2$. More precisely, an {\em approximation} is an infinite connected, convex substructure $A$ of $P_i$ ($i = 1$ or 2) which 
is a finite non-empty union of maximal chains, all isomorphic to $1 + {\cal Z} + 1$ and copies of $\Delta$, such that if $z_1 < z_2$ in such a copy, 
then $z_1$ and $z_2$ are consecutive in $P_i$, and if $z$ is an upward ramification point of $A$ then $z$ is a minimal point of a copy of $\Delta$ in 
$A$, if it is a downward ramification point of $A$ then it is a maximal point of a copy of $\Delta$ in $A$, together with a homomorphism onto 
$1 + {\cal Z} + 1$ which restricts to an isomorphism on each maximal chain. To use this class of approximations to verify uniqueness, it suffices to 
show how to extend any isomorphism $\theta: A \to B$ where $A \subseteq P_1$ and $B \subseteq P_2$ are approximations to include any given point of 
$P_1$ in its domain, and to include any given point of $P_2$ in its range. The two are proved in the same way, so we just do the first. 

Let $A$ and $B$ lie in this class of approximations, with $A \cong B$. Let $A'$ be some extension of $A$ in the class. Then we must show that there 
is $B' \supseteq B$ such that the isomorphism extends. We suppose $A \subset A' \subseteq P_1$ and $B \subseteq P_2$. There are various cases to 
consider, and $A'$ can be formed from $A$ by repeating these cases. 

\noindent{\bf Case 1}: $A'$ is formed from $A$ by adding a copy of $\Delta$, by which is meant that there are $z_1 < z_2$ in $Z$ and $C \cong \Delta$ 
such that $A' = A \cup C$ and $A \cap C = \{z_1, z_2\}$. 

It follows from the definition of `approximation' that $z_1$ and $z_2$ are consecutive in $P_1$. By property (P5), this copy of $\Delta$ is its only 
copy having $(z_1, z_2)$ as an edge, and it again follows from the definition of `approximation' that $z_1$ is not an upward ramification point of 
$A$ and $z_2$ is not a downward ramification point of $A$. Since $\theta$ is an isomorphism, these properties carry across to $B$ in $P_2$, and by 
property (P5) there is a unique $D \cong \Delta$ such that $B \cap D = \{\theta z_1, \theta z_2\}$, and by 1-arc-transitivity of $\Delta$ there is an 
isomorphism from $C$ to $D$ which takes $\{z_1, z_2\}$ to $\{\theta z_1, \theta z_2\}$. Since $A \cap C = \{z_1, z_2\}$ and $B \cap D = \{\theta z_1, 
\theta z_2\}$, the union $\theta'$ of $\theta$ and this isomorphism is a bijection from $A'$ to $B'$. We see that it is an isomorphism thus. Suppose 
that $z_3 < z_4$ in $A'$. If $z_3$ and $z_4$ either both lie in $A$, or both lie in $C$, then $\theta' z_3 < \theta' z_4$ follows from the fact that 
$\theta'$ is an isomorphism on each of these sets, so we suppose that $z_3 \in C \setminus A$ and $z_4 \in A \setminus C$ (the proof if 
$z_4 \in C \setminus A$ and $z_3 \in A \setminus C$ being similar). Now there are two paths from $z_3$ to $z_4$, one `directly upwards' by the fact 
that $z_3 < z_4$, and the other using the connectedness of $A$ and $C$, via $z_1$ or $z_2$. Concatenating these and deleting duplicated sections gives 
a non-alternating cycle, contrary to (P2), unless $z_3 < z_2 < z_4$. But as $\theta'$ is an isomorphism on each of $A$ and $C$, 
$\theta' z_3 < \theta' z_2$ and $\theta' z_2 < \theta' z_4$, from which we deduce $\theta' z_3 < \theta' z_4$. Similarly, 
$\theta' z_3 < \theta' z_4 \Rightarrow z_3 < z_4$.

\noindent{\bf Case 2}: $A'$ is formed from $A$ by adding a maximal point $y$ of $P_1$ above some point $z$ which is maximal in $A$ and in some copy of 
$\Delta$ in $A$, and a chain of points in between. 

For this we observe that $\theta z$ is maximal in $B$, but not in $P_2$, so by appealing to property (P6) we can find a point 
$y_2 \in P_2 \setminus B$ which is maximal in $P_2$ above $\theta z$ and use this to extend $\theta$ to $\theta'$ to include $y$ in its domain. Now 
$[z,y] \cong [\theta z, y_2]$ as one sees using the levels function $l$, and $\theta'$ is obtained as the union of this isomorphism with $\theta$. To
see that this is an isomorphism, first note that by maximality of $z$ in $A$ and $\theta z$ in $B$, 
$(z,y] \cap A = (\theta z, y_2] \cap B = \varnothing$, and therefore $\theta'$ is well defined and a bijection. We check that $\theta'$ is an 
isomorphism thus. Suppose $z_3 < z_4$ in $A'$. If both $z_3$ and $z_4$ lie in $A$, or they both lie in $[z,y]$, then $\theta' z_3 < \theta' z_4$ 
follows from the fact that $\theta'$ is an isomorphism on both these sets, so we may suppose otherwise. By convexity of $A$ the only case we need 
consider is $z_3 \in A \setminus [z,y]$ and $z_4 \in [z,y] \setminus A$. We claim that this implies $z_3 < z < z_4$. Indeed, otherwise we could form a 
cycle $C$ in $A'$ by taking $[z_3,z_4]$ together with a path in $A$ from $z_3$ to $z$ and deleting duplicated sections. Now the cycle $C$ need not 
contain the vertex $z_4$, but we must have $C \cap (z,y] \neq \varnothing$, since we are assuming $z_3$ is not below $z$. But $C$ is a cycle some of 
whose points belong to $A$ and ramify in $A$. By definition of approximation this forces $C \subseteq A$ contradicting 
$C \cap (z,y] \neq \varnothing$. This completes the proof that $z_3 < z < z_4$ and it follows that $A'$ is precisely the poset given by taking the 
transitive closure of $A$ and $(z,y]$. We conclude that $\theta'z_3 < \theta'z < \theta'z_4$. Similarly 
$\theta'z_3 < \theta'z_4 \Rightarrow z_3 < z_4$ completing the proof that $\theta'$ is an isomorphism. 

\noindent{\bf Case 3}: $A'$ is formed from $A$ by adding a minimal point of $P_1$ below some point which is minimal in some copy of $\Delta$ in $A$, 
and a chain of points in between. This is essentially the same as Case 2 so is omitted.

\noindent{\bf Case 4}: $A'$ is formed from $A$ by adding a maximal point $y$ of $P_1$ above some point $z_2$ which is maximal in $A$ and in some copy 
of $\Delta$ in $A$, and a chain of points in between, and a minimal point $x$ of $P_1$ below some point $z_1 < z_2$ which is minimal in some copy of 
$\Delta$ in $A$ (which may be the same as the first one mentioned). 

This case is handled by combining the two previous cases.

To see that these suffice, let $x \in A' \setminus A$. As $A'$ is connected, there is a path in $A'$ from $x$ to a point of $A$. We may add the 
`turning points' of the path one at a time. In other words, we may assume that $x$ is comparable with some member of $A$, suppose $x > y \in A$ for 
instance, where $y$ is maximal in $A$ below $x$. As $A$ is convex, $x$ is not below any member of $A$, and we choose $t$ to be maximal in $P_1$ such 
that $t \ge x$. If $y$ is maximal in $A$, then we use Case 2 to extend to include all points of the interval $[y, t]$ in the domain of $\theta$, since 
by definition of $A$ an `approximation', $y$ is maximal in a copy of $\Delta$ in $A$. We check that $A \cup [y,t]$ is an approximation, provided that 
$y$ is also greater than some minimal member of $P_1$ in $A$, and the only clause requiring verification is convexity. Let $u, v \in A \cup [y,t]$ and 
$u < z < v$, and suppose that $u \in A$ and $v \in [y,t]$. There is a direct path from $u$ to $v$, and also one got by appealing to connectedness of 
each of $A$ and $[y,t]$, and as in Case 2, it follows that $u \le y \le v$. Since also $u < z < v$, by diamond-freeness (which follows from (P7)) it 
follows that $y$ and $z$ are comparable, and hence that $z \in A \cup [y,t]$. If $y$ is not greater than any minimal member $P_1$ in $A$, then as $A$ 
is a finite union of maximal chains and copies of $\Delta$, it follows that there is some minimal member $z' \le y$ of $A$ lying in a copy of 
$\Delta$, and we choose a minimal member $t' < z'$ of $P_1$, and add the interval $[t', z']$ as well as $[y, t]$ (which is Case 4).

If however $y$ is not maximal in $A$, then it becomes a ramification point of $A \cup [y, t]$, and so to fulfil the definition of `approximation' we 
first have to extend to make $y$ the lower point of a copy of $\Delta$ using Case 1. More specifically, let $y < z_1, z_2$ where these are distinct 
and consecutive with $y$, and such that $z_1 \in A$ and $z_2 \le x$ (and hence $z_2 \not \in A$). Note that such points exist even in the 
${\mathbb Q}.2$ case, since by (P1) the upper point of any pair does not ramify upwards in $P_1$. By (P5), there is a unique copy $C$ of $\Delta$ in 
$P_1$ containing $\{y, z_1\}$. We have to see that $A \cup C$ is an approximation, $A \cap C = \{y, z_1\}$, and that $z_2 \in C$ with $z_2$ maximal in 
$A \cup C$. It will follow that we can perform an extension as in Case 1 and add $z_2$ to the approximation.

By (P4), there is a copy of $\Delta$ containing all members of $y^\uparrow$ consecutive with $y$. This must contain $\{y, z_1\}$ and so by uniqueness 
of $C$, must equal $C$. It must also contain $z_2$, and therefore $z_2 \in C$. If $z_2$ is not maximal in $A \cup C$, then as it is certainly maximal 
in $C$, there must be some member of $A$ greater than $z_2$, contrary to convexity of $A$, since $z_2 \not \in A$. Next suppose that 
$u \in (A \cap C) \setminus \{y, z_1\}$. Then as $A$ and $C$ are both connected, there are paths in each of them from $y$ to $u$. If these are 
unequal, then this gives rise to a cycle in $P_1$. By (P3) this lies in a copy of $\Delta$, which by (P5) must be $C$. In either case there is 
therefore a path in $A \cap C$ from $y$ to $u$. Assume that such $u$ is chosen on this path at least distance from $y$. If $u > y$ then as $A$ is an 
approximation, and $y < z_1, u$ with $z_1 \neq u$, there is a copy of $\Delta$ containing $z_1$ and $u$, and contained in $A$. But this copy can again 
only be $C$, and so $C \subseteq A$, which gives $z_2 \in A$, and a contradiction. If however $u \not > y$ then we must have $u < z_1$, and $y, u < 
z_1$ with $y \neq u$, so again using $A$ an approximation, $C \subseteq A$. It remains to show that $A \cup C$ is an approximation. As above (in the 
$A \cup [y,t]$ case) the only part requiring verification is convexity, and this is proved by the same method as in that case.

The proof is concluded by a standard back-and-forth argument. Enumerate the members of each of $P_1$ and $P_2$. Start with an isomorphism between 
single maximal chains of $P_1$ and $P_2$ (using property (P7)), and extend in countably many stages, at the $n$th stage ensuring that the $n$th points 
in the enumerations of $P_1$ and $P_2$ lie in the domain and range respectively.  

The fact that $M$ is 3-$CS$-homogeneous is accomplished by an adaptation of the same back-and-forth method (using the local 2-arc-transitivity of 
$\Delta$). Note that this illustrates that we would not want to require an automorphism to preserve the levels function, since this will be in general 
violated at the first step in trying to map one $Y$-configuration to another. The facts that $M$ is connected and countable are immediate, and both 
$\cal Z$ and $\Delta$ may be recovered from $M$ from the order-type of intervals, and the reachability digraph of $M$ respectively.

We remark that properties (P2) and (P3) are both used in the given proof, which we hope makes things clearer, though one of them would suffice, since 
each can be seen to be easily derivable from the other as follows. Clearly (P3) $\Rightarrow$ (P2) since $\Delta$ can only embed alternating 
cycles. Conversely, by (P2) any cycle must be alternating. Let it be $z_0 < z_1 > z_2 < z_3 \ldots < z_{2n-1} > z_0$ say. By (P4), for each $i$ there 
is a copy $\Delta_i$ of $\Delta$ containing $z_{i-1}$, $z_i$, and $z_{i+1}$ (where the subscripts are taken modulo $2n$). By (P5), all the $\Delta_i$ 
are equal since $\Delta_i$ and $\Delta_{i+1}$ share a common edge, so we have one copy of $\Delta$ containing the whole of the cycle.    $\Box$

\vspace{.1in}

Observe that the construction given in Theorem~\ref{3.16} generalizes the universal highly arc-transitive digraphs construction $DL(\Delta)$ since in the case $\mathcal{Z} = \mathbb{Z}$ the poset $P$ constructed in the theorem, viewed as a digraph, is precisely the digraph $DL(\Delta)$. Combining Theorem~\ref{3.16} with the observation after the proof of Proposition \ref{3.6} that there are $2^{\aleph_0}$ pairwise non-isomorphic countable  $DM$-complete locally $2$-arc-transitive bipartite graphs, we obtain the following. 

\begin{corollary}\label{3.17}
Let $\cal Z$ be a linear order equal to $\mathbb{Z}^{\alpha}$ 
or $\mathbb{Q}.\mathbb{Z}^{\alpha}$ for some countable ordinal $\alpha \ge 1$, or $\mathbb{Q}.2$. 
Then there are $2^{\aleph_0}$ pairwise non-isomorphic countable connected $3$-$CS$-homogeneous $2$-level partial orders $M$ such that the interval of $M$ has order-type $1 + \mathcal{Z} + 1$.
\end{corollary} 

We remark that this corollary justifies the claim made in the abstract that we have, in Theorem \ref{3.16}, described a new family of countably 
infinite locally 2-arc-transitive graphs, namely one for each of the 1-transitive linear orders of the stated form, each containing $2^{\aleph_0}$ 
members.

\begin{corollary}\label{3.18}
If $P$ is the partially ordered set whose existence is given by Theorem \ref{3.16}, there is a homomorphism from $P$ onto $1 + {\cal Z} + 1$. 
\end{corollary}
\noindent{\bf Proof}:  The homomorphism was built into the construction, for the existence proof, and so it exists in that case, but by uniqueness, it 
must exist just on the basis of the properties (P1)-(P8). $\Box$

\vspace{.1in}

Observe that the method of construction outlined above can be used to obtain any of the 3-$CS$-homogeneous $CFPO$s from Warren's classification 
\cite{Warren} for those values of the interval by taking $\Delta = ALT$, the alternating line digraph (i.e. the unique countable connected bipartite graph all of whose vertices have degree equal to two), as input. 
 
\begin{remark}\label{3.19}
{\rm Let $M$ be a countable $3$-$CS$-homogeneous $2$-level connected poset such that $I(M)$ is an infinite chain. In Lemma~\ref{3.9} we determined up 
to isomorphism all the possibilities for $I(M)$. The only cases in Lemma \ref{3.9} which have not so far been discussed are the dense ones, namely 
$1 + \mathbb{Q} + 1$ and $1 + {\mathbb Q}_2 + 1$. We believe that it should be possible to construct many examples $M$ where for every pair 
$x,y \in M^+$ where $x$ is minimal in $M^+$, $y$ is maximal, and $x<y$ we have $[x,y] \cong 1 + \mathbb{Q} + 1$ or $1 + {\mathbb Q}_2 + 1$. The details 
in these cases will however be considerably more involved than the ones so far described, and we have not yet succeeded in working these out. In any
case, the flavour is radically altered, and this is because we cannot any more sensibly regard the intermediate structures as digraphs in any 
meaningful way, because of their density, so for this reason alone, their study belongs elsewhere. 

Now if we attempt to embed $\Delta$ in $Z$, then it cannot be as a convex subset, so instead we should try to embed `extended' versions of $\Delta$, 
where each edge of $\Delta$ is replaced by a copy of the rational interval $[0, 1]$ (2-coloured in the case of ${\mathbb Q}_2$). The difficulty comes 
about because there is now no clear reason for embedding this extended copy of $\Delta$ between any particular two levels of $Z$, so we may be 
obliged to embed it in all possible ways. This forces up the required ramification order, and makes it hard to arrange uniqueness.

A similar but possibly less acute problem arises if we wish to consider embedding $\Delta$ even in the $1 + {\mathbb Z} + 1$ case, but now not 
necessarily on consecutive levels (where the `extended' version of $\Delta$ could have its maximal chains paths of length 2 instead of 1, for 
instance). Some results about this situation are given in \cite{Chen}.} \end{remark}

\section{Completions in general, further examples, and conclusions}

We saw in Proposition~\ref{3.6} that if $I = I(M)$ is a finite chain then $I$ can have at most three elements, and that the resulting examples are incidence graphs of semilinear spaces satisfying a certain transitivity condition. In this section we make some observations concerning countable locally 2-arc-transitive bipartite graphs whose intervals are not chains, concentrating mainly on the case where the interval is finite. 

In Theorem~21 of \cite{Gray2} we proved that the only connected 2-arc-transitive bipartite graph with diamond intervals is the $4 \times 4$ 
complement of perfect matching. In particular, there is no infinite 2-arc-transitive bipartite graph whose completion has diamond intervals.
This leads naturally to the analysis of the case where the interval $I(M)$ has maximal chains of length 3. 

\begin{lemma}\label{4.1}
Let $M$ be a $2$-level partial order. Then any member $x$ of $M^D \setminus M$ is above at least $2$ members of $X = \min(M)$ and is below at least 
$2$ members of $Y = \max(M)$.
\end{lemma}
\noindent{\bf Proof}: Viewing $x$ as an ideal in $M$, $x \subseteq X$, $x, x^\uparrow \neq \emptyset$, and $x^{\uparrow\downarrow} = x$. Since 
$x \neq \emptyset$, $x$ has some member of $X$ below it but if this was unique, $x$ would lie in $X$. Since $x^\uparrow \neq \emptyset$ we can choose 
$v \in Y$, $v \ge x$. If $x^\uparrow = \{v\}$, then $v \in x^{\uparrow\downarrow} = x$, so $x \in Y \subseteq M$. Hence there must be some other 
member of $Y$ above $x$. $\Box$

\vspace{.1in}

Let $M$ be a 2-level 2- and 3-$CS$-transitive countable connected poset. We begin with some general observations. 

The following lemma tells us that meets of pairs of points all lie on the `second top level' of $M^D$. We write $b >> a$ to mean that $b > a$ and 
$(a,b) = \varnothing$ (i.e. there are no points in between). 

\begin{lemma}\label{4.2}
Suppose that the interval $I$ of $M^D$ is finite, and let $y = b \wedge c$ where $b, c \in \max(M)$. Then $b>>y$ and $c>>y$. 
Conversely if $b$ and $c$ are distinct members of $\max(M)$ and there is $y$ such that $b>>y$ and $c>>y$ then $b \wedge c = y$. 
\end{lemma}
\noindent{\bf Proof}:
Suppose for a contradiction that $(y,b) \neq \emptyset$. Then since $I$ is finite we can choose $x \in (y, b)$ so that $b >> x$. Let $d \in \max(M)$ 
with $d>x$ but $d \neq b$. Let $a \in \min(M)$ with $a<y$. Then by 3-$CS$-transitivity there is an automorphism $\theta$ taking $\{a, b, c\}$ to 
$\{a, b, d\}$. The natural extension of $\theta$ to $M^D$ fixes $a$ and sends $y = b \wedge c$ to $b \wedge d = x$, so it acts non-trivially on a 
finite chain, which is impossible.

Conversely, suppose that $b>>y$ and $c>>y$ where $b \neq c$. Then $y \le b \wedge c < b, c$, and so the only possibility is $y = b \wedge c$. $\Box$

\begin{lemma}\label{4.3}
Let $M$ be a $2$-level $2$- and $3$-$CS$-transitive partial order such that for some $a < b$ in $M$ and $x$ in $M^D$, $a < x < b$ and $(a, x)$ and 
$(x, b)$ are empty (in $M^D$). Then any maximal chain of $M^D$ has length $3$.
\end{lemma}
\noindent{\bf Proof}:
Since all intervals are isomorphic, it suffices to prove this for $[a, b]^{M^D}$. If not, there are $y < z$ in $(a, b)^{M^D}$. By Lemma \ref{4.1} let 
$c, d$ be different from $b$ and such that $x < c$, $z < d$. By 3-$CS$-transitivity, there is an automorphism $\theta$ of $M$ taking $\{a, b, c\}$ to 
$\{a, b, d\}$. This fixes $a$ and takes $\{b, c\}$ to $\{b, d\}$. Now $\theta(x) < b, d$ so $z, \theta(x) < b, d$ and hence by completeness of $M^D$ 
there is $t \in M^D$ such that $z, \theta(x) \le t \le b, d$. As $b \neq d$, $t < b, d$. As $x$ is adjacent to $b$, $\theta(x)$ is adjacent to $b$ or 
$d$, from which it follows that $\theta(x) = t$. Hence $z \le \theta(x)$. But as $a$ is adjacent to $x$, it is also adjacent to $\theta(x)$, contrary 
to $a < y < z$.  $\Box$

\subsection{Local finiteness result}

Let $M$ be a poset.  
We say that $M$ is \emph{locally finite} if (i) all intervals $[a,b]_M$ are finite and (ii) the Hasse graph $\Gamma(M)$ associated with $M$ is a 
locally finite graph. (Note that in Theorem \ref{3.15} we referred to `local finiteness', but there it was for digraphs, where the meaning is well 
known---all vertices have finite in- and out-degree.) Let $C_n$ denote the $n$-element chain (that is $\{1, \ldots, n \}$ with the usual ordering). 

\begin{prop}\label{4.4}
Let $M$ be a $2$- and $3$-CS-transitive $2$-level poset and let $I$ be the interval of $M^D$, and suppose that $I \not\cong C_2$ and $I \not\cong C_3$. 
Then $M$ is locally finite if and only if $I$ is finite. 
\end{prop}
\noindent{\bf Proof}:
One direction holds even without any transitivity assumption. For if $M$ is locally finite, then the intervals of $M^D$ are certainly finite (just by 
the definition of $DM$-completion as the family of non-empty bounded above ideals of $M$).

For the converse, suppose that $I$ is finite. Let $a \in \min(M)$ and $b \in \max(M)$ with $a < b$. We shall first show that for any $x \in (a, b)$ 
there is $x' \in (a, b)$ incomparable with $x$. Now by assumption that $[a, b] \not \cong C_2, C_3$, there is $y \neq x$ in $(a, b)$. If $y \parallel 
x$ then we may let $x' = y$. Otherwise $y < x$ or $y > x$, suppose the former. Then there is $c \in \max(M)$ above $y$ but not above $x$, and we let 
$x' = b \wedge c$. Clearly $a < x' < b$, and $x \not \le x'$ since $x \not \le c$ and $x \not > x'$ by Lemma \ref{4.2}.

We can now show that the graph $\Gamma(M)$ is locally finite, and hence, $M$ is locally finite. If not, then there is $a \in \min(M)$ with infinite 
degree (or $b \in \max(M)$ with infinite degree, which is done by a similar argument). Let $\{b_i : i \in \mathbb{N} \}$ be an infinite subset of the 
neighbourhood $\Gamma(a)$. Then since $[a, b_0]$ is finite it follows that the set $A = \{b_0 \wedge b_j : j > 0 \}$ is finite, so there is an 
infinite subset $E$ of $\Gamma(a)$ and $z \in [a,b_0]$ such that $b_0 \wedge b_j = z$ for all $b_j \in E$, and it follows easily that $e \wedge f = z$ 
for all $e, f \in E$, $e \neq f$. By the previous paragraph there is $q \in (a,b_0)$ with $q \parallel z$. Let $c \in \max(M)$ with $c>q$ but 
$c \not \ge z$. Then $c \not\in E$ since $c \not \ge z$. Since $[a, c]$ is finite there are $e, f \in E$, $e \neq f$, such that 
$c \wedge e = c \wedge f = z'$. But $e \wedge f = z$, so since $z' < e,f$ it follows that $z \geq z'$. Since $c \not \ge z$ but $c > z'$ we cannot 
have $z = z'$, so
\[
e \wedge f = z > z' = c \wedge e
\]
and this contradicts Lemma \ref{4.2}. Hence $M$ is locally finite. 
$\Box$

\vspace{.1in}

Note that the assumption that $I$ is not $C_2$ or $C_3$ is necessary. This is seen by inspecting the sporadic examples in Warren's classification (see 
\cite{Warren}). Also, the results from \cite{Gray2} show that the result above is far from being true if we replace 3-$CS$-transitivity by 
2-$CS$-transitivity. 

We think the following is also true, but at the moment it is just a conjecture. 

\begin{conjecture}\label{4.5}
Let $M$ be a connected 3-$CS$-transitive 2-level poset and let $I$ be the interval of $M^D$, and suppose that $I \not\cong C_2$ and $I \not\cong C_3$. 
Then $M$ is finite if and only if $I$ is finite. 
\end{conjecture}

\subsection{Some small cases in detail}

By a \emph{$k$-diamond} we mean a poset with elements $a,b, c_1, \ldots, c_k$ with $a < c_i < b$ (for all $i$) and $a<b$, the only non-trivial 
relations. In particular a 2-diamond is just a diamond, and a $k$-diamond is a diamond with an antichain of $k$ elements on the middle level. 

We now describe those locally finite examples whose completion has intervals isomorphic to 3-diamonds. For this we have to make a technical assumption 
on the diameter (which it may be possible to remove). The argument of Proposition~\ref{4.4} can be used to show that the upward and downward 
ramification orders are each no larger than 3. 

\vspace{.1in}

\noindent{\bf (i) Ramification order 2}

By the `incidence graph of the block design arising from the complement of the Fano plane' we mean a bipartite graph with one part corresponding to 
the lines of the Fano plane, the other part corresponding to the points of the Fano plane, and a line-vertex adjacent to a point-vertex if and only if 
the point does \emph{not} lie on the line in the Fano plane. 

\begin{prop}\label{4.5}
Let $M$ be a locally finite connected $1$- and $2$-arc-transitive bipartite graph such that the intervals of the completion $M^D$ are isomorphic to a 
$3$-diamond. If the upward and downward ramification orders of the midlevel points of $M^D$ are both equal to $2$, and $M$ has diameter at most $3$, 
then $M$ is the incidence graph of the block design arising from the complement of the Fano plane.   
\end{prop}
\noindent{\bf Proof}: We just give an outline.

Let $X = \min(M)$, $Y = \max(M)$ and $Z = M^D \setminus M$, so $Z$ is the set of midlevel points of $M^D$. 

The first step is to describe the structure of $x^\uparrow$ (which now denotes the set of points of $M^D$ which are $\ge x$) for each lower point $x$, 
and show that $x$ has 6 points above it in $Z$ and 4 points above it in $Y$, using a counting argument (and the hypothesis on the intervals). In fact, 
the argument allows us to completely determine the structure of $x^\uparrow$ for any $x \in X$. To `draw a picture' of $x^\uparrow$ one should carry out the following steps:
\begin{enumerate}
\item Start with the minimal point $x$.
\item Then draw the 6 points immediately above $x$ in the middle level. 
\item Draw 4 points above the 6 points (which will be the 4 maximal points of $x^\uparrow$).
\item Finally, for each distinct pair of the 4 points on the top level define the meet of this pair to be one of the 6 points in the middle level 
(with distinct pairs giving distinct meets), and draw in lines for these relations (this will give each of the 6 midlevel points upward ramification 
order 2 in $x^\uparrow$). 
\end{enumerate}
If one draws this picture then it is easy to see that $x^\uparrow$ has the following property:
\begin{quote}
Every point $u$ in $x^\uparrow \cap Z$ (i.e. midlevel point) has a unique mate $u' \in x^\uparrow \cap Z$ such that $u$ and $u'$ do \emph{not} have a 
common upper bound in $x^\uparrow$. 
\end{quote}

In the conclusion of the proof, one considers $\mathcal{B} = \{ \Gamma(y) : y \in Y \}$ (where $\Gamma(y)$ denotes the set of all neighbours in $M$ of 
the vertex $y$) which is a set of 4-subsets of $X$, and shows that it is a $(v,b,r,k,\lambda)$-design with $k = r = 4$ and $\lambda = 2$. From 
$\lambda(v-1) = r(k-1)$ and $bk = vr$ it follows that $v = b = 7$, so $|X| = |Y| = 7$ and $M$ is a $7 \times 7$ bipartite graph. The bipartite 
complement $\bar{M}$ of the graph $M$ (so $\bar{M}$ has the same vertex set as $M$ and the same bipartition, and $\{x,y\}$ is an edge in $\bar{M}$ if 
and only if it is not an edge in $M$) is then seen to be a projective plane, and since the Fano plane is the unique projective plane of order 2 (see 
\cite{Veblen}) it follows that $\bar{M}$ is isomorphic to the incidence graph of the Fano plane, completing the proof of the theorem. $\Box$


\vspace{.1in}

\noindent{\bf (ii) Ramification order 3}


There is an example in this case also. First we describe the example, and then we prove it is the unique such example. 

Let $V(n,q)$ be a vector space of dimension $n$ ($n \geq 2$) over the field $\mathbb{F}_q$ with $q=p^m$ elements ($p$ prime). Let $X$ be the set of 
1-dimensional subspaces of $V$ and let $Y$ be the set of $(n-1)$-dimensional subspaces of $V$. Let $M=M(n,q)$ be the bipartite graph with vertex set 
$X \cup Y$ and $x \sim y$ if and only if $x \subseteq y$. Then for all $n$ and $q=p^m$ the graph $M(n,q)$ is a $2$-arc-transitive bipartite graph 
whose completion $M^D$ has $n-1$ levels, since it is easily seen that the members of $M^D$ are the intermediate subspaces of $V(n,q)$. 

In particular the completion of $M(4,q)$ has $3$ levels, and its intervals are isomorphic to a $k$-diamond where $k = (q^3-q) / (q^2-q) = q+1$ which is the number of $2$-dimensional subspaces of a fixed $3$-dimensional subspace of $V(n,q)$ containing some fixed $1$-dimensional subspace. Also, the upward and downward ramification order of the points on the middle level of $M(4,q)$ are both equal to $(q^4 - q^2)/(q^3-q^2)=q+1$. 

So, for any prime power $q$ there is an example whose interval is a $(q+1)$-diamond and whose ramification order is $q+1$. In particular, when $q=2$ we have a $3$-diamond with ramification orders $3$. 

\begin{lemma}\label{4.6}
Let $M$ be a connected $2$-arc-transitive bipartite graph such that the intervals of the completion $M^D$ are isomorphic to a $k$-diamond. Suppose 
that the upward and downward ramification orders of the midlevel points both equal $r$. If $k=r$ then for all $a \in \min{M}$ the number of midlevel 
points in $a^{\uparrow}$ and the number of maximal points in  $a^{\uparrow}$ are both equal to $r(r-1)+1$, and any two points in the middle level of  
$a^{\uparrow}$ have a common upper bound in $a^{\uparrow}$. 
\end{lemma}
\noindent{\bf Proof}:
Let $X = \min(M)$, $Y = \max(M)$ and $Z = M^D \setminus (X \cup Y)$. Fix $a \in X$ and let $U = a^\uparrow \cap Z$ and $V = a^\uparrow \cap Y$. 
First we shall prove that $|V| = r(r-1)+1$. 

Fix some $b \in V$. For every $v \in V \setminus \{b\}$, since $v$ and $b$ are bounded below by $a$ it follows that $b \wedge v$ exists and 
$b \wedge v \in [a,b] \cap U$. Thus every $v \in V$ lies above one of the $r$ points in $[a,b] \cap U$. Along with the assumption that midlevel points 
have upward ramification order $r$, this implies that $|V| \leq r(r-1) + 1$. On the other hand, every point $v$ of $V \setminus \{b\}$ lies above a 
unique member of $[a,b] \cap U$, as follows from the $DM$-completeness of $M^D$, since if $v > u_1, u_2$ where $u_1$ and $u_2$ are distinct members of 
$U$, $v \ge u_1 \vee u_2 = b$. This shows that $|V| \geq r(r-1)+1$, completing the proof that $|V| = r(r-1)+1$. 

Next, since the upward ramification order of each midlevel point is $r$, it follows that every point of $U = a^{\uparrow} \cap Z$ can be written as $x 
\wedge y$ where $x,y \in Y$ in exactly ${r \choose 2} = r(r-1)/2$ different ways. Hence, taking meets of all possible pairs in $V$ we deduce that 
\[
|U| = 
{r(r-1)+1 \choose 2} / {r \choose 2} = [(r(r-1)+1)(r(r-1)) / 2][2 / r(r-1)]
= r(r-1)+1.
\]
This completes the proof of the first part of the lemma.

For the second part, view $U \cup V$ as a bipartite graph. Since $|U| = |V|$ and all members of $U$ have degree $r$, so do all members of $V$. Hence 
the number of vertices of $U$ at distance 2 from a given $u \in U$ is equal to $r(r-1)$, since $u$ is adjacent to exactly $r$ points of $V$, 
each of which is adjacent to $r-1$ other points of $U$ (and all these points are distinct by $DM$-completeness). As $U$ has size $r(r-1)+1$ and $u$ 
was arbitrary, this completes the proof of the lemma. $\Box$

\begin{prop}\label{4.7}
Let $M$ be a connected $2$-arc-transitive bipartite graph such that the intervals of the completion $M^D$ are isomorphic to a $3$-diamond. If the 
upward and downward ramification orders of the midlevel points of $M^D$ are both equal to $3$ then $M \cong M(4,2)$. 
\end{prop}
\noindent{\bf Proof}:
Let $x$ and $y$ be in the lower level and suppose that $d_M(x,y)=4$. Let $x,a,z,b,y$ be a path of length 4 from $x$ to $y$. Then $x \vee z$ and $z 
\vee y$ are in $z^{\uparrow}$ so by Lemma~\ref{4.6} they have a common upper bound in $z^{\uparrow}$. This upper bound also serves as a common 
upper bound for the pair $x,y$ contradicting the fact that $d_M(x,y)=4$. We conclude that any two lower points of $M$ have a common upper bound in $M$.

Now we want to recover the whole structure of $M$. Let $X = \min(M)$ and $Y = \max(M)$. Let $\mathcal{B} = \{ \Gamma(y) : y \in Y \}$ which is a set 
of 7-subsets of $X$ (using Lemma \ref{4.6}, and since $r = k = 3$, so $r(r-1) + 1 = 7$). Since for all $x_1, x_2 \in X$, $d(x_1,x_2)=2$ and $M$ is 
2-arc-transitive it follows that the automorphism group $G$ of $M$ acts 2-transitively on $Y$, and so is doubly transitive on the blocks 
$\mathcal{B}$, and any two blocks $B_i, B_j$ intersect. So there is a number $N$ such that for any two blocks $B_i, B_j$ we have $|B_i \cap B_j| = N$. 
It follows quickly from the assumptions on $M$ that $N = 3$. So, $X$ is a finite set, and $\mathcal{B}$ is a set of 7-element subsets of $X$ such that 
each point belongs to 7 blocks. Moreover any two blocks intersect in a set of size 3, and any pair of distinct points belongs to exactly 3 blocks. 
This gives us a $(15,7,3)$-design, since from $\lambda(v-1) = r(k-1)$ it follows that $v = b = 15$, so $|X| = |Y| = 15$ and $M$ is a $15 \times 15$ 
bipartite graph. 

But Kantor \cite{Kantor} classified all the 2-transitive symmetric designs and from his classification a 2-transitive $(15,7,3)$-design must be a 
projective space. The only projective space with 15 points is $M(4,2)$. $\Box$

\vspace{.1in}

Currently we do not know if this generalizes to give the following. 

\

\noindent \textbf{Question.} Let $q$ be a prime power. Let $M$ be a connected $2$-arc-transitive bipartite graph such that the intervals of the completion $M^D$ are isomorphic to a $(q+1)$-diamond. If the upward and downward ramification orders of the midlevel points of $M^D$ are both equal to $q+1$ then does it follow that $M \cong M(4,q)$?  

\subsection{Additional examples}

Given our general approach outlined in Subection~\ref{sec_outline}, one other natural line of investigation is to begin with a family of bipartite graphs that are known to be locally 2-arc-transitive and investigate the resulting Dedekind--MacNeille completions. In this section we list a few well-known such families and make some observations about the structure of their Dedekind--MacNeille completions.

\vspace{.1in}

\noindent{\bf (i) Generalized cubes}. 

We view a {\em generalized cube} as a bipartite graph. We write it as $C_n$, and this exists for all positive integers $n$, and also in 
the infinite case (even uncountably infinite). For ease we start with the finite case, and take $C_n$ to consist of all binary sequences of length 
$n$. We take $X$ and $Y$ to consist of the sequences with an even number (odd number, respectively) of 1s, and we join $\sigma$ and $\tau$ if they 
differ on exactly one entry. In the infinite case, we start with an infinite cardinal $\kappa$, and this time work with the set of all functions 
$\sigma$ from $\kappa$ into $\{0, 1\}$ such that $\sigma^{-1}\{1\}$ is finite; the definitions of $X$ and $Y$ and adjacency are as in the finite case. 
We remark that it is possible to construct the finite $C_n$ inductively, and there is an obvious way to embed $C_n$ in $C_{n+1}$, and the infinite 
cases may be formed as suitable limits. The main point of the examples however, is that each $C_n$ is 1- and 2-arc transitive, and the interval of its 
completion is an $n$-diamond (whether $n$ is finite or infinite). 

To verify the claimed transitivity of $C_n$, we may use automorphisms of $C_n$ of two types. The first are ones which are induced by permutations of 
$\{0, 1, \ldots, n-1\}$, and the second are ones which swap 0 and 1 entries on a fixed set of entries of even size; the details are omitted.

Finally to check that the intervals $I$ in the completion of $C_n$ are $(n-1)$-diamonds, the main things to be done is to show that the members of  
$C_n^D \setminus C_n$ lie on one level, and they are precisely the 2-element subsets of $X$ whose elements differ in exactly 2 places. Given
$\sigma < \tau$ in $C_n$, which differ in the $i$th position, there are therefore exactly $n-1$ members of $C_n^D \setminus C_n$ between them, namely 
subsets of $X$ of the form $\{\sigma, \sigma'\}$ where $\sigma'$ differs from $\sigma$ at positions $i$ and $j$ for some $j \neq i$. Thus the interval 
$I$ is an $(n-1)$-diamond. (This applies in the infinite case too.) 

\vspace{.1in}

\noindent{\bf (ii) The complement of a perfect matching}

Let $X$ be any non-empty set, and let $Y$ be a disjoint copy of $X$. We may regard $X \cup Y$ as a bipartite graph where $x$ is joined to $y$ if $x 
\neq y$. Since a `perfect matching' is the similar bipartite graph in which $x$ is joined to $y \Leftrightarrow$ $x = y$, we have here the `complement 
of a perfect matching'. The fact that this is locally 2-arc-transitive follows from the stronger property that it is homogeneous (any isomorphism 
between finite substructures which respects the parts extends to an automorphism). One can check that the completion is isomorphic to the power set of 
$X$ (with $\emptyset$ and $X$ removed). 

Note that this example applies whether $X$ is finite or infinite. The case where $|X| = 4$ was studied in \cite{Gray2} (Theorem 25).

\vspace{.1in}

\noindent{\bf (iii) The generic  bipartite graph}

This is the bipartite graph on $X \cup Y$ characterized as follows. $X$ and $Y$ are disjoint countably infinite sets, and for each finite disjoint 
$U_1, U_2 \subseteq X$ there is $y \in Y$ joined to all of $U_1$ and none of $U_2$, and similarly if $U_1, U_2$ are finite disjoint subsets of $Y$, 
there is a corresponding element of $X$. To characterize precisely what $(X \cup Y)^D$ is seems quite hard, but we know many of its members. For 
instance, all finite subsets of $X$ are ideals, as follows easily from genericity, and similarly for cofinite subsets. This is not all however, as it 
is easy to construct infinite and coinfinite ideals, using a diagonal argument. To describe the structure of the partial ordering of all such sets 
seems very complicated.

\vspace{.1in}

\noindent{\bf (iv) A 2-arc-transitive bipartite graph derived from subspaces of a vector space}

Let $V$ be an $n$-dimensional vector space over a field $F$, where $n \geq 3$. Let $X$ and $Y$ be the families of subspaces of $V$ of dimension, 
codimension 1, respectively, and let $x<y$ if $x \not\subseteq y$, where $x \in X$, $y \in Y$. Then one can show that $X \cup Y$ is a locally 
2-arc-transitive bipartite graph, using some straightforward linear algebra.

Now we examine the intervals in some special cases, so we have to see what the points of $(X \cup Y)^D$ are. The first remark that for any 
${\bf x}_1 \neq {\bf x}_2$ in $X$, $\{{\bf x}_1, {\bf x}_2\} \in (X \cup Y)^D$. 

We now treat separately the cases where $F$ has 2 members or more. If $F = F_2$, the field with 2 elements, then we can show that the members of 
$(X \cup Y)^D$ not in $Y$ are precisely the subsets of $X$ of the form `the set of sums of an odd number of elements of 
$\{{\bf x}_1, {\bf x}_2, \ldots, {\bf x}_r\}$' where ${\bf x}_1, {\bf x}_2, \ldots, {\bf x}_r$ are linearly independent. This shows that for $F_2$, 
$(X \cup Y)^D$ has exactly $n$ levels, and we can in principle work out how many points there are on each level. For instance, if $n = 3$ there are 7, 
21, 7 points on the three levels (being the example discussed in Lemma \ref{4.6}) and if $n = 4$ there are 15, 105, 105, 15. The intervals in 
these cases are a projective line and a projective plane over $F_2$ (the projective line, when viewed as a lattice, being a 3-diamond).

If $|F| > 2$, the members of $(X \cup Y)^D$ are more complicated to describe. In particular, if ${\bf x}_1, {\bf x}_2, {\bf x}_3$ are distinct, then 
one shows that $\{{\bf x}_1, {\bf x}_2, {\bf x}_3\} \in (X \cup Y)^D$ (by different proofs depending on whether the three vectors are linearly 
dependent or not), but we have not determined what the general form a a member of the completion takes in this example.. 

\subsection{Conclusions}

We conclude by listing the main types of 3-$CS$-transitive diamond-free but not cycle-free partial orders that we have constructed or discussed in this 
paper. 

The largest and most significant family is as given by Theorem \ref{3.16}. For each countable ordinal $\alpha \ge 1$, and for each countable connected 
$DM$-complete locally 2-arc-transitive bipartite graph $\Delta$ there are corresponding 2-level 3-$CS$-homogeneous partial orders, whose maximal chains 
have order-type $1 + {\mathbb Z}^\alpha + 1$, $1 + {\mathbb Q}.{\mathbb Z}^\alpha + 1$, and ${\mathbb Q}.2$ and whose reachability digraph is $\Delta$.

As a general remark, it is shown that the finite chain 3-$CS$-homogeneous partial orders which are Dedekind--MacNeille complete are precisely given by the incidence graphs of semilinear spaces which fulfil the natural corresponding transitivity condition.

The other examples, given in this final section, serve as a preliminary investigation into the case where the interval is not a chain, and give some information about the intervals that arise from certain classical families of $2$-arc-transitive bipartite graphs such as 
generalized cubes, the complement of a perfect matching, the generic bipartite graph, and various examples derived from subspaces of 
vector spaces.

Authors' addresses

  \medskip
Robert Gray, \\
School of Mathematics, University of East Anglia\\
Norwich Research Park, Norwich NR4 7TJ UK

 \medskip

 \texttt{Robert.D.Gray@uea.ac.uk} \\
    
    John K. Truss, \\
    Department of Pure Mathematics, \\
    University of Leeds, Leeds LS2 9JT, UK \\
    
\texttt{pmtjkt@leeds.ac.uk}

\end{document}